\documentclass[11pt]{article}
\usepackage{epsfig}

\def\be{\begin{equation}}
\def\ee{\end{equation}}
\def\bea{\begin{eqnarray}}
\def\eea{\end{eqnarray}}
\def\bes{\begin{eqnarray*}}
\def\ees{\end{eqnarray*}}

\def\nn{\nonumber}
\def\lb{\label}
\def\bs{\setminus}

\def\vs{{\varsigma}}

\def\R{{\bf R}}
\def\C{{\bf C}}
\def\Z{{\bf Z}}
\def\K{{\bf K}}
\def\N{{\bf N}}
\def\U{{\bf U}}

\def\Q{{\bf Q}}

\def\aa{{\alpha}}
\def\bb{{\beta}}
\def\ga{{\gamma}}

\def\th{{\theta}}

\def\om{{\omega}}
\def\Om{{\Omega}}
\def\ep{{\epsilon}}
\def\lm{{\lambda}}
\def\Lm{{\Lambda}}

\def\sg{{\sigma}}
\def\dm{{\diamond}}

\def\vf{{\varphi}}

\def\<{{\langle}}
\def\>{{\rangle}}

\def\K{{\cal K}}
\def\P{{\cal P}}

\def\M{{\cal M}}
\def\Nn{{\cal N}}

\def\rank{{\rm rank}}

\def\Sp{{\rm Sp}}

\def\mod{{\rm mod}}
\def\ol{\overline}
\def\td#1{\tilde{#1}}

\def\hb{\vrule height0.18cm width0.14cm $\,$}

\title{Multiplicity and stability of closed geodesics on positively curved Finsler $4$-spheres}

\author{Huagui Duan\thanks{Partially supported by National Key R\&D Program of China (Grant No. 2020YFA0713300), NNSFC (Nos. 12271268 and 11790271), and the Fundamental Research Funds for the Central Universities. E-mail: duanhg@nankai.edu.cn.}, \quad
Dong Xie\thanks{Partially supported by NNSFC (No. 12271268). E-mail: 2120170062@mail.nankai.edu.cn.} \\ \\
School of Mathematical Sciences and LPMC, Nankai University,\\ Tianjin 300071, P. R. China\\}

\begin{document}

\maketitle

\begin{abstract}
{\it In this paper, we prove that for every Finsler $4$-dimensional sphere $(S^4,F)$ with reversibility $\lambda$ and flag curvature $K$ satisfying $\frac{25}{9}\left(\frac{\lambda}{1+\lambda}\right)^2<K\le 1$ with $\lambda<\frac{3}{2}$, either there exist at least four prime closed geodesics, or there exist exactly three prime non-hyperbolic closed geodesics and at least two of them are irrationally elliptic.}
\end{abstract}

{\bf Key words}: Positively curved, closed geodesic, irrationally elliptic, Finsler metric, spheres.

{\bf 2020 Mathematics Subject Classification}: 53C22, 58E05, 58E10.

\renewcommand{\theequation}{\thesection.\arabic{equation}}
\renewcommand{\thefigure}{\thesection.\arabic{figure}}

\setcounter{figure}{0}
\setcounter{equation}{0}
\section{Introduction and main result}

Let $(M,F)$ be a Finsler manifold. A closed curve on $(M,F)$ is a closed geodesic if it is
locally the shortest path connecting any two nearby points on this
curve. As usual, a closed geodesic $c:S^1=\R/\Z\to M$ is {\it prime}
if it is not a multiple covering (i.e., iteration) of any other
closed geodesics. Here the $m$-th iteration $c^m$ of $c$ is defined
by $c^m(t)=c(mt)$. The inverse curve $c^{-1}$ of $c$ is defined by
$c^{-1}(t)=c(1-t)$ for $t\in \R$.  Note that unlike Riemannian manifold,
the inverse curve $c^{-1}$ of a closed geodesic $c$
on an irreversible Finsler manifold need not be a geodesic.
We call two prime closed geodesics
$c$ and $d$ {\it distinct} if there is no $\th\in (0,1)$ such that
$c(t)=d(t+\th)$ for all $t\in\R$.
On a reversible Finsler (or Riemannian) manifold, two closed geodesics
$c$ and $d$ are called { \it geometrically distinct} if $
c(S^1)\neq d(S^1)$, i.e., they have different image sets in $M$.
We shall omit the word {\it distinct} when we talk about more than one prime closed geodesic.

For a closed geodesic $c$ on an $(n+1)$-dimensional manifold $M$, denote by $P_c$ the linearized Poincar\'{e} map of $c$, which is a symplectic matrix, i.e., $P_c\in\Sp(2n)$. we define the {\it elliptic height } $e(P_c)$ of $P_c$ to be the total algebraic multiplicity of all eigenvalues of $P_c$ on the unit circle $\U=\{z\in\C|\; |z|=1\}$ in the complex plane $\C$. Since $P_c$ is symplectic, $e(P_c)$ is even and $0\le e(P_c)\le 2n$.
A closed geodesic $c$ is called {\it elliptic} if $e(P_c)=2n$, i.e., all the eigenvalues of $P_c$ locate on $\U$; {\it irrationally elliptic} if, in the homotopy component $\Om^0(P_c)$ of $P_c$ (cf. Section 2 below for the definition), $P_c$ can be connected to the $\dm$-product of $n$ rotation matrices $R(\th_i)$ with $\th_i$ being irrational multiple of $\pi$ for $1\le i\le n$; {\it hyperbolic} if $e(P_c)=0$, i.e., all the
eigenvalues of $P_c$ locate away from $\U$; {\it non-degenerate} if $1$ is not an eigenvalue of $P_c$. A Finsler metric $F$
is called {\it bumpy} if all the closed geodesics on $(M,F)$ are non-degenerate.

There is a famous conjecture in Riemannian geometry which claims the existence of infinitely many
closed geodesics on any compact Riemannian manifold. This conjecture
has been proved for many cases, but not yet for compact rank one symmetric
spaces except for $S^2$. The results of Franks in \cite{Fra} and Bangert in \cite{Ban} imply that this conjecture is true for any Riemannian 2-sphere (cf. \cite{Hin2} and \cite{Hin3}). However for a Finsler manifold, the above conjecture does not hold due to the Katok's examples.
It was quite surprising when Katok in \cite{Kat} found some irreversible Finsler metrics on spheres with only finitely
many closed geodesics and all of them are non-degenerate and irrationally elliptic (cf. \cite{Zil}).

Based on Katok's examples, Anosov in \cite{Ano} proposed the following conjecture (cf. \cite{Lon4})
\bea \Nn(S^n,F)\ge 2\left[\frac{n+1}{2}\right] \quad \mbox{for any Finsler metric $F$ on}\ S^n, \lb{1.1}\eea
where, denote by $\Nn(M,F)$ the number of the distinct closed geodesics on $(M,F)$, and $[a]=\max\{k\in\Z\,|\,k\le a\}$. In 2005, Bangert and Long in \cite{BaL} proved this conjecture for any Finsler $2$-dimensional sphere $(S^2, F)$. Since then, the index iteration theory of closed geodesics (cf. \cite{Bot} and \cite{Lon3}) has been applied to study the closed geodesic problem on Finsler manifolds.

When $n\ge 3$, the above conjectures in the Riemannian or Finsler case is still widely open in full of generality. About the multiplicity and stability problem of closed geodesics, two classes of typical conditions, including the positively curved condition and the non-degenerate (or bumpy) condition, have been used widely.

In \cite{Rad3}, Rademacher has introduced the reversibility $\lambda=\lambda(M,F)$ of a compact Finsler manifold defined by
\bea \lambda=\max\{F(-X)\ |\ X\in TM,\ F(X)=1\}\ge 1.\nn\eea
Then Rademacher in \cite{Rad4} has obtained some results about the multiplicity and stability of closed geodesics. For example, let $F$ be a Finsler metric on $S^{n}$ with reversibility $\lm$ and flag curvature $K$ satisfying $\left(\frac{\lm}{1+\lm}\right)^2<K\le 1$, then there exist at least $n/2-1$ closed geodesics with length $<2n\pi$. If $\frac{9\lm^2}{4(1+\lm)^2}<K\le 1$ with $\lm<2$, then there exists a closed geodesic of elliptic-parabolic, i.e., its linearized Poincar\'{e} map split into $2$-dimensional rotations and a part whose eigenvalues are $\pm 1$. These results are some generalization of those in \cite{BTZ1} and \cite{BTZ2} in the Riemannian case.

Recently, Wang in \cite{Wan} proved the conjecture (\ref{1.1}) for $(S^n,F)$ provided that $F$ is bumpy and its flag curvature $K$ satisfies $\left(\frac{\lambda}{1+\lambda}\right)^2<K\le 1$. Also in \cite{Wan}, Wang showed that for every bumpy Finsler metric $F$ on $S^n$ satisfying $\frac{9\lm^2}{4(1+\lm)^2}<K\le 1$, there exist two prime elliptic closed geodesics provided the number of closed geodesics on $(S^n,F)$ is finite.
As some further generalization, Duan, Long and Wang in \cite{DLW} obtained the optimal lower bound of the number of distinct closed geodesics on a compact simply-connected Finsler manifold $(M,F)$ if $F$ is bumpy and some much weak index conditions or positively curved conditions are satisfied.

The first author in \cite{Dua1} and \cite{Dua2} proved that for every Finsler $(S^{n},F)$ for $n\ge 3$ with reversibility $\lm$ and flag curvature $K$ satisfying $\left(\frac{\lm}{1+\lm}\right)^2<K\le 1$, either there exist infinitely many prime closed geodesics, or there exist exactly three prime closed geodesics and at least two of them are elliptic. In fact, the multiplicity and stability problem on high dimensional manifolds without the assumption of bumpy metrics is much difficult.

In this paper, we further consider the positively curved Finsler $4$-dimensional sphere $(S^4,F)$ without the bumpy assumption, and obtain the following new progress about the multiplicity and stability of closed geodesics on $(S^4,F)$.

\medskip

{\bf Theorem 1.1.} {\it For every Finsler metric $F$ on a $4$-dimensional sphere $S^4$ with reversibility $\lm$ and flag curvature $K$ satisfying $\frac{25}{9}\left(\frac{\lambda}{1+\lambda}\right)^2<K\le 1$ with $\lambda<\frac{3}{2}$, either there exist at least four prime closed geodesics, or there exist exactly three prime non-hyperbolic closed geodesics and at least two of them are irrationally elliptic.}

\medskip

First, under the assumption of positively curved condition $\frac{25}{9}\left(\frac{\lambda}{1+\lambda}\right)^2<K\le 1$ with $\lambda<\frac{3}{2}$, Theorem 1 and Theorem 4 in \cite{Rad3} established the lower bound of the length of closed geodesics, which in turn gives the lower bound of $i(c^m)$ and mean index $\hat{i}(c)$ for any prime closed geodesic $c$ on such $S^4$ (cf. Lemma 3.1 below). Second, we shall make full use of the enhanced common index jump theorem established in \cite{DLW}, which generalized the common index jump theorem in \cite{LoZ}, to obtain some crucial precise estimates of $i(c^m)$ and $\nu(c^m)$ (cf. Section 3.1 below). Note that Theorem 1.1 in \cite{Dua1} showed the existence of three prime closed geodesics on $(S^4,F)$ with reversibility $\lm$ and flag curvature $K$ satisfying $\left(\frac{\lm}{1+\lm}\right)^2<K\le 1$. So, in order to prove our above Theorem 1.1, we assume the existence of exactly three prime closed geodesics on such $(S^4,F)$ (cf. the assumption (TCG) below). Finally, together with the Morse theory, under (TCG) we will carefully analyze the local and global information of these three prime closed geodesics and their iterates to complete the proof of Theorem 1.1 in Section 3.2.

In addition, under the assumption (TCG), in Theorem 4.2 in Section 4, we obtain more precise information about the third prime closed geodesic except that it is non-hyperbolic. These information will be greatly helpful to completely solve the conjecture (\ref{1.1}) on the positively cured Finsler $(S^4,F)$ in the future.

\medskip

{\bf Remark 1.2.} (i) The irrationally ellipticity of closed geodesics is a kind of very important stability, in fact, it specially implies the non-degeneracy of closed geodesics. It is conjectured that all closed geodesics are irrationally elliptic if the number of prime closed geodesics on Finsler sphere $(S^n, F)$ with $n\ge 2$ is finite (cf. Conjecture of \cite{DLW}). In this direction, there are several well-known results. For example, Long and Wang in \cite{LW} proved that on every Finsler $S^2$ with only finitely many closed geodesics, there exist at least two irrationally elliptic one. Duan and Liu in \cite{DL16} showed that if there exist exactly three prime closed geodesics on every Finsler $S^3$ satisfying $\frac{9}{4}\left(\frac{\lambda}{1+\lambda}\right)^2<K\le 1$ with $\lambda<2$, then two
of them are irrationally elliptic.

(ii) This paper continues to consider the case of $n=4$ of $S^n$. Compared to previous works, and in particular to \cite{Dua1} and \cite{Dua2}, the first substantial progress in this paper is to prove the existence of two prime irrationally elliptic closed geodesics on the above positively curved Finsler $S^4$ if there exist exactly three prime closed geodesics on it. Except for this, the second progress in this paper is to obtain five possible characterizations about the Morse index and the linearized Poincar\'{e} map of the third closed geodesic if there exist exactly three prime closed geodesics on such $S^4$ (cf. Theorem 4.2 below). These results is very difficult due to the high dimension of $n\ge 4$. In fact, in the case of $n=4$, there possibly exist three $2\times 2$ rotation matrices in the decomposition of $P_c$ (cf. (\ref{2.15}) of Theorem 2.6 below), which brings high degeneracy of closed geodesics so that some methods and technical tools in \cite{Dua1}, \cite{Dua2} and \cite{DL16} are not enough to deal with these difficulties. In this paper, we make full use of the enhanced common index jump theorem established in \cite{DLW}, together with some new ideas and estimates about indices of closed geodesics.

(iii) On one hand, note that the curvature pinching condition in Theorem 1.1 and Theorem 4.2 in this paper is a little restrictive than the corresponding assumption in \cite{Dua1} and \cite{Dua2}, this is because that we need good index estimates in Lemma 3.1 and Lemma 3.2 to count the contributions of iterated closed geodesics to the non-zero Morse-type numbers, and the current pinching condition is a sufficient one for these estimates. On the other hand, we do not consider the case $n\ge 5$ because of more higher degeneracy. This maybe need some new ideas and different tools to deal with more difficulties.

\medskip

In this paper, let $\N$, $\N_0$, $\Z$, $\Q$, $\R$, and $\C$ denote the sets of natural integers, non-negative integers, integers,
rational numbers, real numbers, and complex numbers respectively. We use only singular homology modules with $\Q$-coefficients.
For an $S^1$-space $X$, we denote by $\overline{X}$ the quotient space $X/S^1$. We define the functions
\bea
E(a)=\min\{k\in\Z\,|\,k\ge a\},\quad \varphi(a)=E(a)-[a],\quad \{a\}=a-[a].  \lb{1.2}
\eea
Especially, $\varphi(a)=0$ if $ a\in\Z\,$, and $\varphi(a)=1$ if $
a\notin\Z\,$.

\setcounter{equation}{0}
\section{Morse theory and Morse indices of closed geodesics}

\subsection{Morse theory for closed geodesics} 

Let $M=(M,F)$ be a compact Finsler manifold, the space
$\Lambda=\Lambda M$ of $H^1$-maps $\gamma:S^1\rightarrow M$ has a
natural structure of Riemannian Hilbert manifolds on which the
group $S^1=\R/\Z$ acts continuously by isometries. This action is defined by
$(s\cdot\gamma)(t)=\gamma(t+s)$ for all $\gamma\in\Lm$ and $s,
t\in S^1$. For any $\gamma\in\Lambda$, the energy functional is
defined by
\be E(\gamma)=\frac{1}{2}\int_{S^1}F(\gamma(t),\dot{\gamma}(t))^2dt.
\lb{2.1}\ee
It is $C^{1,1}$ and invariant under the $S^1$-action. The
critical points of $E$ of positive energies are precisely the closed geodesics
$\gamma:S^1\to M$. The index form of the functional $E$ is well
defined along any closed geodesic $c$ on $M$, which we denote by
$E''(c)$. As usual, we denote by $i(c)$ and
$\nu(c)$ the Morse index and nullity of $E$ at $c$. In the
following, we denote by
\be \Lm^\kappa=\{d\in \Lm\;|\;E(d)\le\kappa\},\quad \Lm^{\kappa-}=\{d\in \Lm\;|\; E(d)<\kappa\},
  \quad \forall \kappa\ge 0. \nn\ee
For a closed geodesic $c$ we set $ \Lm(c)=\{\ga\in\Lm\mid E(\ga)<E(c)\}$.

Recall that respectively the mean index $\hat{i}(c)$ and the $S^1$-critical modules of $c^m$ are defined by
\be \hat{i}(c)=\lim_{m\rightarrow\infty}\frac{i(c^m)}{m}, \quad \overline{C}_*(E,c^m)
   = H_*\left((\Lm(c^m)\cup S^1\cdot c^m)/S^1,\Lm(c^m)/S^1; \Q\right).\lb{2.3}\ee

We call a closed geodesic satisfying the isolation condition, if
the following holds:

\medskip

{\bf (Iso)  For all $m\in\N$ the orbit $S^1\cdot c^m$ is an
isolated critical orbit of $E$. }

\medskip

Note that if the number of prime closed geodesics on a Finsler manifold
is finite, then all the closed geodesics satisfy (Iso).

If $c$ has multiplicity $m$, then the subgroup $\Z_m=\{\frac{n}{m}\mid 0\leq n<m\}$
of $S^1$ acts on $\overline{C}_*(E,c)$. As studied in p.59 of \cite{Rad2},
for all $m\in\N$, let
$H_{\ast}(X,A)^{\pm\Z_m}
   = \{[\xi]\in H_{\ast}(X,A)\,|\,T_{\ast}[\xi]=\pm [\xi]\}$,
where $T$ is a generator of the $\Z_m$-action.
On $S^1$-critical modules of $c^m$, the following lemma holds:

\medskip

{\bf Lemma 2.1.} (cf. Satz 6.11 of \cite{Rad2} and \cite{BaL})
{\it Suppose $c$ is
a prime closed geodesic on a Finsler manifold $M$ satisfying (Iso). Then
there exist $U_{c^m}^-$ and $N_{c^m}$, the so-called local negative
disk and the local characteristic manifold at $c^m$ respectively,
such that $\nu(c^m)=\dim N_{c^m}$ and
\bea \overline{C}_q( E,c^m)
&\equiv& H_q\left((\Lm(c^m)\cup S^1\cdot c^m)/S^1, \Lm(c^m)/S^1\right)\nn\\
&=& \left(H_{i(c^m)}(U_{c^m}^-\cup\{c^m\},U_{c^m}^-)
    \otimes H_{q-i(c^m)}(N_{c^m}\cup\{c^m\},N_{c^m})\right)^{+\Z_m}, \nn
\eea

(i) When $\nu(c^m)=0$, there holds
$$ \overline{C}_q( E,c^m) = \left\{\matrix{
     \Q, &\quad {\it if}\;\; i(c^m)-i(c)\in 2\Z\;\;{\it and}\;\;
                   q=i(c^m),\;  \cr
     0, &\quad {\it otherwise}, \cr}\right.  $$

(ii) When $\nu(c^m)>0$, there holds
$$ \overline{C}_q( E,c^m)=H_{q-i(c^m)}(N_{c^m}\cup\{c^m\},N_{c^m})^{\ep(c^m)\Z_m}, $$
where $\ep(c^m)=(-1)^{i(c^m)-i(c)}$.}

\medskip

Define
\be  k_j(c^m) \equiv \dim\, H_j( N_{c^m}\cup\{c^m\},N_{c^m}), \quad
     k_j^{\pm 1}(c^m) \equiv \dim\, H_j(N_{c^m}\cup\{c^m\},N_{c^m})^{\pm\Z_m}. \lb{2.4}\ee

Then we have

\medskip

{\bf Lemma 2.2.} (cf. \cite{Rad2}, \cite{LoD}, \cite{Wan})
{\it Let $c$ be a  prime closed geodesic on a Finsler manifold $(M,F)$. Then

(i) For any $m\in\N$, there holds
$k_j(c^m)=0$ for $j\not\in [0,\nu(c^m)]$.

(ii) For any $m\in\N$, $k_0(c^m)+k_{\nu(c^m)}(c^m)\le 1$ and if $k_0(c^m)+k_{\nu(c^m)}(c^m)=1$ then there holds  $k_j(c^m)=0$ for $j\in (0,\nu(c^m))$.

(iii) For any $m\in\N$, there holds $k_0^{+1}(c^m) = k_0(c^m)$ and $k_0^{-1}(c^m) = 0$.
In particular, if $c^m$ is non-degenerate, there holds $k_0^{+1}(c^m) = k_0(c^m)=1$,
and $k_0^{-1}(c^m) = k_j^{\pm 1}(c^m)=0$ for all $j\neq 0$.

(iv) Suppose for some integer $m=np\ge 2$ with $n$ and $p\in\N$ the nullities
satisfy $\nu(c^m)=\nu(c^n)$. Then there holds $k_j(c^m)=k_j(c^n)$ and
${k}_j^{\pm 1}(c^m)={k}^{\pm 1}_j(c^n)$ for any integer $j$. }

\medskip

Let $(M,F)$ be a compact simply connected Finsler manifold with finitely many closed geodesics. It
is well known that for every prime closed geodesic $c$ on $(M,F)$, there holds either
$\hat{i}(c)>0$ and then $i(c^m)\to +\infty$ as $m\to +\infty$, or $\hat{i}(c)=0$ and
then $i(c^m)=0$ for all $m\in\N$. Denote those prime closed geodesics on $(M,F)$ with
positive mean indices by $\{c_j\}_{1\le j\le k}$.
Rademacher in \cite{Rad} and \cite{Rad2} established a celebrated mean index identity relating all the $c_j$s
with the global homology of $M$ for compact simply connected Finsler manifolds (especially for $S^4$) as follows.

\medskip

{\bf Theorem 2.3.} (Satz 7.9 of \cite{Rad2}, cf. also \cite{DuL}, \cite{LoD} and \cite{Wan})
{\it Assume that there exist finitely many prime closed geodesics on $(S^4,F)$ and
denote prime closed geodesics with positive mean indices by
$\{c_j\}_{1\le j\le k}$ for some $k\in\N$. Then the following identity holds
\be  \sum_{j=1}^k\frac{\hat{\chi}(c_j)}{\hat{i}(c_j)} = -\frac{2}{3}, \lb{2.5}\ee
where
\be \hat{\chi}(c_j)
= \frac{1}{n(c_j)}\sum_{1\le m\le n(c_j)}
     \chi(c_j^m)=\frac{1}{n(c_j)}\sum_{1\le m\le n(c_j) \atop 0\le l\le 2(n-1)}(-1)^{i(c_j^m)+l}k_l^{\ep(c_j^m)}(c_j^m)\in\Q,    \lb{2.6}\ee
and the analytical period
$n(c_j)$ of $c_j$ is defined by (cf. \cite{LoD})
\be n(c_j) = \min\{l\in\N\,|\,\nu(c_j^l)=\max_{m\ge 1}\nu(c_j^m),\;\;
                  i(c_j^{m+l})-i(c_j^{m})\in 2\Z, \;\;\forall\,m\in\N\}. \lb{2.7}\ee}

\medskip

Set $\ol{\Lm}^0=\ol{\Lambda}^0S^4 =\{{\rm constant\;point\;curves\;in\;}S^4\}\cong S^4$.
Let $(X,Y)$ be a space pair such that the Betti numbers $b_i=b_i(X,Y)=\dim H_i(X,Y;\Q)$
are finite for all $i\in \Z$. As usual the {\it Poincar\'e series} of $(X,Y)$ is
defined by the formal power series $P(X, Y)=\sum_{i=0}^{\infty}b_it^i$. We need the
following well known version of results on Betti numbers and the Morse inequality.

\medskip

{\bf Lemma 2.4.} (cf. Theorem 2.4 and Remark 2.5 of \cite{Rad} and \cite{Hin1}, Lemma 2.5 of \cite{DuL})
{\it Let $(S^4,F)$ be a $4$-dimensional Finsler sphere. Then, the Betti numbers are given by
\bea
b_j
&=& \rank H_j(\Lm S^4/S^1,\Lm^0 S^4/S^1;\Q)  \nn\\
&=& \left\{\matrix{
    2,&\quad {\it if}\quad j\in \K\equiv \{3k\,|\,3\le k\in 2\N+1\},  \cr
    1,&\quad {\it if}\quad j\in \{2k+3\,|\,k\in\N_0\}\bs\K,  \cr
    0,&\quad {\it otherwise}. \cr}\right.    \lb{2.8}
\eea}

{\bf Theorem 2.5.} (cf. Theorem I.4.3 of \cite{Cha})
{\it Let $(M,F)$ be a Finsler manifold with finitely many prime closed geodesics, denoted by $\{c_j\}_{1\le j\le k}$. Set
\bea M_q =\sum_{1\le j\le k,\; m\ge 1}\dim{\ol{C}}_q(E, c^m_j),\quad q\in\Z.\nn\eea
Then for every integer $q\ge 0$ there holds }
\bea M_q - M_{q-1} + \cdots +(-1)^{q}M_0
&\ge& {b}_q - {b}_{q-1}+ \cdots + (-1)^{q}{b}_0, \lb{2.9}\\
M_q &\ge& {b}_q. \lb{2.10}\eea

\subsection{Index iteration theory of closed geodesics}

In \cite{Lon1} of 1999, Y. Long established the basic normal form
decomposition of symplectic matrices. Based on this result he
further established the precise iteration formulae of indices of
symplectic paths in \cite{Lon2} of 2000. Note that this index iteration formulae works for Morse indices of iterated closed geodesics (cf. \cite{Liu} and Chap. 12 of \cite{Lon3}). Since every closed geodesic on a sphere must be orientable. Then by Theorem 1.1 of \cite{Liu}, the initial Morse index of a closed geodesic on a Finsler $S^4$ coincides with the index of a corresponding symplectic path.

As in \cite{Lon2}, denote by
\bea
N_1(\lm, b) &=& \left(\matrix{\lm & b\cr
                                0 & \lm\cr}\right), \qquad {\rm for\;}\lm=\pm 1, \; b\in\R, \lb{2.11}\\
D(\lm) &=& \left(\matrix{\lm & 0\cr
                      0 & \lm^{-1}\cr}\right), \qquad {\rm for\;}\lm\in\R\bs\{0, \pm 1\}, \lb{2.12}\\
R(\th) &=& \left(\matrix{\cos\th & -\sin\th \cr
                           \sin\th & \cos\th\cr}\right), \qquad {\rm for\;}\th\in (0,\pi)\cup (\pi,2\pi), \lb{2.13}\\
N_2(e^{\th\sqrt{-1}}, B) &=& \left(\matrix{ R(\th) & B \cr
                  0 & R(\th)\cr}\right), \qquad {\rm for\;}\th\in (0,\pi)\cup (\pi,2\pi)\;\; {\rm and}\; \nn\\
        && \qquad B=\left(\matrix{b_1 & b_2\cr
                                  b_3 & b_4\cr}\right)\; {\rm with}\; b_j\in\R, \;\;
                                         {\rm and}\;\; b_2\not= b_3. \lb{2.14}\eea
Here $N_2(e^{\th\sqrt{-1}}, B)$ is non-trivial if $(b_2-b_3)\sin\theta<0$, and trivial
if $(b_2-b_3)\sin\theta>0$.

As in \cite{Lon2}, the $\diamond$-sum (direct sum) of any two real matrices is defined by
$$ \left(\matrix{A_1 & B_1\cr C_1 & D_1\cr}\right)_{2i\times 2i}\diamond
      \left(\matrix{A_2 & B_2\cr C_2 & D_2\cr}\right)_{2j\times 2j}
=\left(\matrix{A_1 & 0 & B_1 & 0 \cr
                                   0 & A_2 & 0& B_2\cr
                                   C_1 & 0 & D_1 & 0 \cr
                                   0 & C_2 & 0 & D_2}\right). $$

For every $M\in\Sp(2n)$, the homotopy set $\Omega(M)$ of $M$ in $\Sp(2n)$ is defined by
$$ \Om(M)=\{N\in\Sp(2n)\,|\,\sg(N)\cap\U=\sg(M)\cap\U\equiv\Gamma\;\mbox{and}
                    \;\nu_{\om}(N)=\nu_{\om}(M),\ \forall\om\in\Gamma\}, $$
where $\sg(M)$ denotes the spectrum of $M$,
$\nu_{\om}(M)\equiv\dim_{\C}\ker_{\C}(M-\om I)$ for $\om\in\U$.
The component $\Om^0(M)$ of $P$ in $\Sp(2n)$ is defined by
the path connected component of $\Om(M)$ containing $M$.

\medskip

{\bf Theorem 2.6.} (cf. Theorem 7.8 of \cite{Lon1}, Theorems 1.2 and 1.3 of \cite{Lon2}, cf. also
Theorem 1.8.10, Lemma 2.3.5 and Theorem 8.3.1 of \cite{Lon3}) {\it For every $P\in\Sp(2n-2)$, there
exists a continuous path $f\in\Om^0(P)$ such that $f(0)=P$ and
\bea f(1)
&=& N_1(1,1)^{\dm p_-}\,\dm\,I_{2p_0}\,\dm\,N_1(1,-1)^{\dm p_+}
  \dm\,N_1(-1,1)^{\dm q_-}\,\dm\,(-I_{2q_0})\,\dm\,N_1(-1,-1)^{\dm q_+} \nn\\
&&\dm\,N_2(e^{\aa_{1}\sqrt{-1}},A_{1})\,\dm\,\cdots\,\dm\,N_2(e^{\aa_{r_{\ast}}\sqrt{-1}},A_{r_{\ast}})
  \dm\,N_2(e^{\bb_{1}\sqrt{-1}},B_{1})\,\dm\,\cdots\,\dm\,N_2(e^{\bb_{r_{0}}\sqrt{-1}},B_{r_{0}})\nn\\
&&\dm\,R(\th_1)\,\dm\,\cdots\,\dm\,R(\th_{r'})\,\dm\,R(\th_{r'+1})\,\dm\,\cdots\,\dm\,R(\th_r)\dm\,H(2)^{\dm h},\lb{2.15}\eea
where $\frac{\th_{j}}{2\pi}\in\Q\cap(0,1)\bs \{\frac{1}{2}\} $ for $1\le j\le r'$ and
$\frac{\th_{j}}{2\pi}\notin\Q\cap(0,1)$ for $r'+1\le j\le r$; $N_2(e^{\aa_{j}\sqrt{-1}},A_{j})$'s
are nontrivial and $N_2(e^{\bb_{j}\sqrt{-1}},B_{j})$'s are trivial, and non-negative integers
$p_-, p_0, p_+,q_-, q_0, q_+,r,r_\ast,r_0,h$ satisfy the equality
\be p_- + p_0 + p_+ + q_- + q_0 + q_+ + r + 2r_{\ast} + 2r_0 + h = n-1. \lb{2.16}\ee

Let $\ga\in\P_{\tau}(2n-2)=\{\ga\in C([0,\tau],\Sp(2n-2))\,|\,\ga(0)=I\}$. we extend $\ga(t)$ to $t \in [0,m\tau]$ for every $m \in \N$ by
\bea
\ga^{m}(t)=\ga(t-j\tau)\ga(\tau)^j \qquad\forall j\tau\le t\le (j+1)\tau\ \mbox{and}\ j=0,1,\cdots,m-1.\lb{2.17}
\eea
Denote the basic normal form
decomposition of $P\equiv \ga(\tau)$ by (\ref{2.15}). Then we have
\bea i(\ga^m)
&=& m(i(\ga)+p_-+p_0-r ) + 2\sum_{j=1}^r{E}\left(\frac{m\th_j}{2\pi}\right) - r   \nn\\
&&  - p_- - p_0 - {{1+(-1)^m}\over 2}(q_0+q_+)
              + 2\sum_{j=1}^{r_{\ast}}\vf\left(\frac{m\aa_j}{2\pi}\right) - 2r_{\ast}, \lb{2.18}\\
\nu(\ga^m)
 &=& \nu(\ga) + {{1+(-1)^m}\over 2}(q_-+2q_0+q_+) + 2\vs(m,\ga(\tau)),    \lb{2.19}\eea
where we denote by
\be \vs(m,\ga(\tau)) = r - \sum_{j=1}^r\vf(\frac{m\th_j}{2\pi})
             + r_{\ast} - \sum_{j=1}^{r_{\ast}}\vf(\frac{m\aa_j}{2\pi})
             + r_0 - \sum_{j=1}^{r_0}\vf(\frac{m\bb_j}{2\pi}).    \lb{2.20}\ee}

\medskip

Let
\be
\M\equiv\{N_1(1,1); \;\;N_1(-1,a_2),\,a_2=\pm1;\;\;R(\th), \,\th\in[0,2\pi);\,H(-2)\}.
\ee

By Theorems 8.1.4-8.1.7 and 8.2.1-8.2.4 of \cite{Lon3}, we have

\medskip

{\bf Proposition 2.7.} {\it Every path $\ga\in\P_{\tau}(2)$ with end matrix homotopic
to some matrix in $\M$ has odd index $i(\ga)$. Every path $\xi\in\P_{\tau}(2)$
with end matrix homotopic to $N_1(1,-1)$ or $H(2)$, and every path $\eta\in\P_{\tau}(4)$
with end matrix homotopic to $N_2(\om,B)$ has even indices $i(\xi)$ and $i(\eta)$.}

\medskip

The common index jump theorem (cf. Theorem 4.3 of \cite{LoZ}) for symplectic paths has become
one of the main tools in studying the multiplicity and stability of periodic orbits in Hamiltonian
and symplectic dynamics. Recently, the following enhanced common index jump theorem has been
obtained by Duan, Long and Wang in \cite{DLW}.

\medskip

{\bf Theorem 2.8.} (cf. Theorem 3.5 of \cite{DLW}) {\it Let
	$\gamma_k\in\mathcal{P}_{\tau_k}(2n)$ for $k=1,\cdots,q$ be a finite collection of symplectic paths.
	Let $M_k=\ga_k(\tau_k)$. We extend $\ga_k$ to $[0,+\infty)$ by (\ref{2.17}) inductively. Suppose
	\be  \hat{i}(\ga_k,1) > 0, \qquad \forall\ k=1,\cdots,q.  \lb{2.21}\ee
	Then for any fixed integer $\bar{m}\in \N$, there exist infinitely many $(q+1)$-tuples
	$(N, m_1,\cdots,m_q) \in \N^{q+1}$ such that for all $1\le k\le q$ and $1\le m\le \bar{m}$, there holds
	\bea
	\nu(\ga_k,2m_k-m) &=& \nu(\ga_k,2m_k+m) = \nu(\ga_k, m),   \lb{2.22}\\
	i(\ga_k,2m_k+m) &=& 2N+i(\ga_k,m),                         \lb{2.23}\\
	i(\ga_k,2m_k-m) &=& 2N-i(\ga_k,m)-2(S^+_{M_k}(1)+Q_k(m)),  \lb{2.24}\\
	i(\ga_k, 2m_k)&=& 2N -(S^+_{M_k}(1)+C(M_k)-2\Delta_k),     \lb{2.25}\eea
	where $S_{M_k}^\pm(\om)$ is the splitting number of $M_k$ at $\om$ (cf. Definition 9.1.4 of \cite{Lon3}) and
	\bea
	&&C(M_k)=\sum\limits_{0<\theta<2\pi}S^-_{M_k}(e^{\sqrt{-1}\theta}),\  \Delta_k = \sum_{0<\{m_k\th/\pi\}<\delta}S^-_{M_k}(e^{\sqrt{-1}\th}),\nn\\
	&&Q_k(m) = \sum_{\theta\in(0,2\pi), e^{\sqrt{-1}\th}\in\sg(M_k),\atop
		\{\frac{m_k\th}{\pi}\}= \{\frac{m\th}{2\pi}\}=0}
	S^-_{M_k}(e^{\sqrt{-1}\th}). \lb{2.26}
	\eea
	More precisely, by (4.10), (4.40) and (4.41) in \cite{LoZ}, we have
	\bea m_k=\left(\left[\frac{N}{\bar{M}\hat i(\gamma_k, 1)}\right]+\chi_k\right)\bar{M},\quad 1\le k\le q,\lb{2.27}\eea
	where $\chi_k=0$ or $1$ for $1\le k\le q$ and $\frac{\bar{M}\theta}{\pi}\in\Z$
	whenever $e^{\sqrt{-1}\theta}\in\sigma(M_k)$ and $\frac{\theta}{\pi}\in\Q$
	for some $1\le k\le q$.  Furthermore, for any fixed $M_0\in\N$, we may
	further require  $M_0|N$, and for any $\epsilon>0$, we can choose $N$ and $\{\chi_k\}_{1\le k\le q}$ such that
	\bea \left|\left\{\frac{N}{\bar{M}\hat i(\gamma_k, 1)}\right\}-\chi_k\right|<\epsilon,\quad 1\le k\le q.\lb{2.28}\eea}

\medskip

We also have the following properties in the index iteration theory.

\medskip

{\bf Theorem 2.9.} (cf. Theorem 2.2 of \cite{LoZ} or Theorem 10.2.3 of \cite{Lon3}) {\it Let $\ga\in P_{\tau}(2n)$. Then, for any $m\in N$, there holds
$$
\nu(\ga,m)-\frac{e(M)}{2}\le i(\ga,m+1)-i(\ga,m)-i(\ga,1)\le\nu(\ga,1)-\nu(\ga,m+1)+\frac{e(M)}{2},
$$
where $e(M)$ is the elliptic height defined in Section 1.}

\setcounter{figure}{0}
\setcounter{equation}{0}
\section{Some index estimates and proof of Theorem 1.1}
\subsection{Some index estimates for closed geodesics}

Firstly we make the following assumption

{\bf (FCG)} {\it Suppose that there exist only finitely many prime closed geodesics $\{c_k\}_{k=1}^q$ on $(S^4,F)$ with reversibility $\lambda$ and flag curvature $K$ satisfying $\frac{25}{9}\left(\frac{\lm}{1 + \lm}\right)^2 < K \le 1$ with $\lambda<\frac{3}{2}$.}

\medskip

For any $1\le k\le q$, we rewrite (\ref{2.15}) as follows
\bea f_k(1)
&=& N_1(1,1)^{\dm p_{k,-}}\,\dm\,I_{2p_{k,0}}\,\dm\,N_1(1,-1)^{\dm p_{k,+}}\nn\\
  &&\dm\,N_1(-1,1)^{\dm q_{k,-}}\,\dm\,(-I_{2q_{k,0}})\,\dm\,N_1(-1,-1)^{\dm q_{k,+}} \nn\\
&& \dm\,R(\th_{k,1})\,\dm\,\cdots\,\dm\,R(\th_{k,r_{k,1}})\,\dm\,R(\td{\th}_{k,1})\,\dm\,\cdots\,\dm\,R(\td{\th}_{k,r_{k,2}})\nn\\
&& \dm\,N_2(e^{\sqrt{-1}\aa_{k,1}},A_{k,1})\,\dm\,\cdots\,\dm\,N_2(e^{\sqrt{-1}\aa_{k,r_{k,3}}},A_{k,r_{k,3}})\nn\\
&& \dm\,N_2(e^{\sqrt{-1}\td{\aa}_{k,1}},\td{A}_{k,1})\,\dm\,\cdots\,\dm\,N_2(e^{\sqrt{-1}\td{\aa}_{k,r_{k,4}}},\td{A}_{k,r_{k,4}})\nn\\
&& \dm\,N_2(e^{\sqrt{-1}\bb_{k,1}},B_{k,1})\,\dm\,\cdots\,\dm\,N_2(e^{\sqrt{-1}\bb_{k,r_{k,5}}},B_{k,r_{k,5}})\nn\\
&& \dm\,N_2(e^{\sqrt{-1}\td{\bb}_{k,1}},\td{B}_{k,1})\,\dm\,\cdots\,\dm\,N_2(e^{\sqrt{-1}\td{\bb}_{k,r_{k,6}}},\td{B}_{k,r_{k,6}})\dm\,H(2)^{\dm h_{k,+}}\dm\,H(-2)^{\dm h_{k,-}},\lb{3.1.0}
\eea
where $\frac{\th_{k,j}}{2\pi}\in\Q\cap(0,1)\bs \{\frac{1}{2}\} $ for $1\le j\le r_{k,1}$, $\frac{\td{\th}_{k,j}}{2\pi}\in(0,1)\bs\Q$ for $1\le j\le r_{k,2}$,
$\frac{\aa_{k,j}}{2\pi}\in\Q\cap(0,1)\bs \{\frac{1}{2}\} $ for $1\le j\le r_{k,3}$, $\frac{\td{\aa}_{k,j}}{2\pi}\in(0,1)\bs\Q$ for $1\le j\le r_{k,4}$,
$\frac{\bb_{k,j}}{2\pi}\in\Q\cap(0,1)\bs \{\frac{1}{2}\} $ for $1\le j\le r_{k,5}$, $\frac{\td{\bb}_{k,j}}{2\pi}\in(0,1)\bs\Q$ for $1\le j\le r_{k,6}$; $N_2(e^{\sqrt{-1}\aa_{k,j}},A_{k,j})$'s and $N_2(e^{\sqrt{-1}\td{\aa}_{k,j}},\td{A}_{k,j})$'s
are nontrivial and $N_2(e^{\sqrt{-1}\bb_{k,j}},B_{k,j})$'s and $N_2(e^{\sqrt{-1}\td{\bb}_{k,j}},\td{B}_{k,j})$'s are trivial, and non-negative integers
$p_{k,-}$, $p_{k,0}$, $p_{k,+}$, $q_{k,-}$, $q_{k,0}$, $q_{k,+}$, $r_{k,1}$, $r_{k,2}$, $r_{k,3}$, $r_{k,4}$, $r_{k,5}$, $r_{k,6}$, $h_k=h_{k,+}+h_{k,-}$ satisfy the equality
\bea
p_{k,-}+p_{k,0}+p_{k,+}+q_{k,-}+q_{k,0}+q_{k,+}+r_{k,1}+r_{k,2} + 2\sum_{j=3}^6 r_{k,j}+h_k= 3. \lb{3.2.0}
\eea

\medskip

{\bf Lemma 3.1.} {\it Under the assumption (FCG), for any prime closed geodesic $c_k$, $1\le k\le q$, there holds
\be
i(c_k^m)\ge 3\left[\frac{5m}{3}\right],\qquad \forall\ m\in\N\lb{3.3.0}
\ee
and
\be
\hat{i}(c_k)>5. \lb{3.4.0}
\ee}

{\bf Proof.} By the assumption (FCG), since the flag curvature $K$ satisfies $\frac{25}{9}\left(\frac{\lm}{1 + \lm}\right)^2 < K \le 1$, we can choose $\frac{25}{9}\left(\frac{\lambda}{\lambda+1}\right)^2<\delta\le K\le 1$. Then by Lemma 2 in \cite{Rad4}, it yields
 $$
 \hat{i}(c_k)\ge 3\sqrt{\delta}\frac{1+\lambda}{\lambda}>5.
 $$

Note that it follows from Theorem 3 of \cite{Rad3} that $L(c_k^m)=mL(c_k)\ge m\pi\frac{1+\lambda}{\lambda}>\frac{5m}{3}\pi/\sqrt{\delta}$ for $m\ge 1$ and $1\le k\le q$. Then it follows from Lemma 3 of \cite{Rad3} that $i(c_k^m)\ge 3[\frac{5m}{3}]$.  \hfill\hb

\medskip

Combining Lemma 3.1 with Theorem 2.9, it follows that
\bea
i(c_k^{m+1})-i(c_k^m)-\nu(c_k^m)\ge i(c_k)-\frac{e(P_{c_k})}{2}\ge 0,\quad\forall\ m\in\N,\ 1\le k\le q.\lb{3.5.0}
\eea
Here the last inequality holds by the fact that $e(P_{c_k})\le 6$ and $i(c_k)\ge 3$.

It follows from (\ref{3.4.0}),  Theorem 4.3 in \cite{LoZ}  and Theorem 2.8 that for any fixed integer $\bar{m}\in \N$, there exist infinitely many $(q+1)$-tuples $(N, m_1,\cdots, m_q)\in\N^{q+1}$ such that for any $1\le k\le q$ and $1\le m\le \bar{m}$, there holds
\bea i(c_k^{2m_k -m})+\nu(c_k^{2m_k-m})&=&
2N-i(c_k^{m})-\left(2S^+_{P_{c_k}}(1)+2Q_{k}(m)-\nu(c_k^{m})\right), \lb{3.6.0}\\
i(c_k^{2m_k})&\ge& 2N-\frac{e(P_{c_k})}{2},\lb{3.7.0}\\
i(c_k^{2m_k})+\nu(c_k^{2m_k})&\le& 2N+\frac{e(P_{c_k})}{2},\lb{3.8.0}\\
i(c_k^{2m_k+m})&=&2N+i(c_k^{m}),\lb{3.9.0}
\eea
where (\ref{3.7.0}) and (\ref{3.8.0}) follow from (4.32) and (4.33) in Theorem 4.3 in \cite{LoZ} respectively.

Note that by List 9.1.12 of \cite{Lon3}, (\ref{2.26}), (\ref{2.19}) and $\nu(c_k)=p_{k,-}+2p_{k,0}+p_{k,+}$, we have
\bea
&&S^+_{P_{c_k}}(1)=p_{k,-}+p_{k,0},\lb{3.10.0}\\
&&C(P_{c_k})=q_{k,0}+q_{k,+}+r_{k,1}+r_{k,2}+2r_{k,3}+2r_{k,4},\lb{3.11.0}\\
&&Q_{k}(m)=\frac{1+(-1)^m}{2}(q_{k,0}+q_{k,+})+(r_{k,1}+r_{k,3}) -\sum_{j=1}^{r_{k,1}}\vf\left(\frac{m\theta_{k,j}}{2\pi}\right)-\sum_{j=1}^{r_{k,3}}\vf\left(\frac{m\aa_{k,j}}{2\pi}\right),\lb{3.12.0}\\
&&\nu(c_k^m)=(p_{k,-}+2p_{k,0}+p_{k,+})+\frac{1+(-1)^m}{2}(q_{k,-}+2q_{k,0}+q_{k,+})+2(r_{k,1}+r_{k,3}+r_{k,5})\nn\\
&&\qquad\qquad -2\left(\sum_{j=1}^{r_{k,1}}\vf\left(\frac{m\th_{k,j}}{2\pi}\right)+\sum_{j=1}^{r_{k,3}}\vf\left(\frac{m\aa_{k,j}}{2\pi}\right)
+\sum_{j=1}^{r_{k,5}}\vf\left(\frac{m\bb_{k,j}}{2\pi}\right)\right).\lb{3.13.0}
\eea

By (\ref{3.10.0}), (\ref{3.12.0}) and (\ref{3.13.0}), we obtain
\bea
2S^+_{P_{c_k}}(1)+2Q_k(m) -\nu(c_k^{m})=p_{k,-}-p_{k,+}-\frac{1+(-1)^m}{2}(q_{k,-}-q_{k,+})
-2r_{k,5}+2\sum_{j=1}^{r_{k,5}}\vf\left(\frac{m\bb_{k,j}}{2\pi}\right),\nn
\eea
which, together with (\ref{3.6.0}), gives
\bea
i(c_k^{2m_k -m})+\nu(c_k^{2m_k-m})&=&
2N-i(c_k^{m})-p_{k,-}+p_{k,+}+\frac{1+(-1)^m}{2}(q_{k,-}-q_{k,+})\nn\\
&&\quad +2r_{k,5}-2\sum_{j=1}^{r_{k,5}}\vf\left(\frac{m\bb_{k,j}}{2\pi}\right),\quad\forall\ 1\le m\le \bar{m}.\lb{3.14.0}
\eea
By (\ref{3.7.0})-(\ref{3.9.0}), (\ref{3.14.0}), (\ref{3.2.0}), (\ref{3.3.0}) and the fact $e(P_{c_k})\le 6$, there holds
\bea
i(c_k^{2m_k-m})+\nu(c_k^{2m_k-m})&\le& 2N+3-3\left[\frac{5m}{3}\right], \quad\forall\ 1\le m\le \bar{m},\lb{3.15.0}\\
2N-3 &\le& 2N-\frac{e(P_{c_k})}{2}\le i(c_k^{2m_k}),\lb{3.16.0}\\
i(c_k^{2m_k})+\nu(c_k^{2m_k})&\le& 2N+\frac{e(P_{c_k})}{2}\le 2N+3,\lb{3.17.0}\\
2N+3\left[\frac{5m}{3}\right]&\le& i(c_k^{2m_k+m}), \qquad\forall\ 1\le m\le \bar{m}.\lb{3.18.0}
\eea

Note that by (\ref{3.5.0}), we have
\bes
&&i(c_k^m)\le i(c_k^{m+1}), \qquad i(c_k^m)+\nu(c_k^m)\le  i(c_k^{m+1})+\nu(c_k^{m+1}),\quad\forall m\in \N,
\ees
which, together with (\ref{3.15.0}) and (\ref{3.18.0}), implies
\bea
&&i(c_k^m)+\nu(c_k^m)\le i(c_k^{2m_k -\bar{m}})+\nu(c_k^{2m_k-\bar{m}})\le 2N+3-3\left[\frac{5\bar{m}}{3}\right], \,\forall\ 1\le m\le 2m_k-\bar{m},\lb{3.19.0}\\
&&2N+3\left[\frac{5\bar{m}}{3}\right]\le i(c_k^{2m_k+\bar{m}})\le i(c_k^{m}),\,\forall\ m\ge 2m_k+\bar{m}.\lb{3.20.0}
\eea

In addition, by (\ref{2.25}), (\ref{3.10.0}), (\ref{3.11.0}) and (\ref{3.13.0}), the precise formulae of $i(c_k^{2m_k})+\nu(c_k^{2m_k})$ can be computed as follows
\bea
i(c_k^{2m_k})+\nu(c_k^{2m_k})
&=&2N+2\Delta_k-(p_{k,-}+p_{k,0}+q_{k,0}+q_{k,+}+r_{k,1}+r_{k,2}+2r_{k,3}+2r_{k,4})\nn\\
&&+p_{k,+}+2p_{k,0}+p_{k,-}+q_{k,+}+2q_{k,0}+q_{k,-}+2r_{k,1}+2r_{k,3}+2r_{k,5}\nn\\
&=&2N+p_{k,0}+p_{k,+}+q_{k,-}+q_{k,0}\nn\\
&&\qquad +r_{k,1}+2r_{k,5}-r_{k,2}-2r_{k,4}+2\Delta_k,\quad k=1,\cdots,q, \lb{3.21.0}
\eea
where
\bea \Delta_k \equiv \sum_{0<\{m_k\th/\pi\}<\delta}S^-_{M_k}(e^{\sqrt{-1}\th})\le r_{k,2}+r_{k,4}. \lb{3.22.0}\eea
by (\ref{2.26}) and List 9.1.12 of \cite{Lon3}.

\medskip

{\bf Lemma 3.2.} {\it Under the assumption (FCG), for $k=1,\cdots,q$ and $1\le m \le \bar{m}$, we have
\bea i(c_k^{2m_k -1})+\nu(c_k^{2m_k -1})&\le& 2N-3,\lb{3.23.0}\\
i(c_k^{2m_k -2})+\nu(c_k^{2m_k -2})&\le& 2N-9.\lb{3.24.0}\eea}

{\bf Proof.} By (\ref{2.18}), we have
\bea \hat i(c_{k})&=&i(c_{k})+p_{k,-}+p_{k,0}-r_{k,1}-r_{k,2}+\sum_{j=1}^{r_{k,1}}
\frac{\theta_{k,j}}{\pi}+\sum_{j=1}^{r_{k,2}}
\frac{\tilde  \theta_{k,j}}{\pi}\nn\\
&<& i(c_{k})+p_{k,-}+ p_{k,0}+r_{k,1}+r_{k,2}.\lb{3.25.0}
\eea
Combining (\ref{3.25.0}) with (\ref{3.4.0}), there holds
\bea i(c_{k})+p_{k,-}+p_{k,0}+r_{k,1}+r_{k,2}\ge 6. \lb{3.26.0}
\eea
Then by (\ref{3.14.0}), (\ref{3.26.0}) and (\ref{3.2.0}), we obtain
\bea
i(c_k^{2m_k -1})+\nu(c_k^{2m_k -1})&=&
2N-i(c_k)-p_{k,-}+p_{k,+} \nn\\
&\le & 2N-6+p_{k,0}+r_{k,1}+r_{k,2}+p_{k,+}\le 2N-3.\lb{3.27.0}
\eea
By (\ref{3.14.0}), (\ref{3.3.0}) and (\ref{3.2.0}), we get
\bea i(c_k^{2m_k -2})+\nu(c_k^{2m_k -2})=
2N-i(c_k^{2})-p_{k,-}+p_{k,+}+q_{k,-}-q_{k,+}\le 2N-9+3=2N-6.\lb{3.28.0}
\eea
Now we assume $ i(c_k^{2m_k -2})+\nu(c_k^{2m_k -2})\ge 2N-8 $, $\forall\ k=1,\cdots,q$, which, together with (\ref{3.28.0}), gives $  i(c_k^{2m_k -2})+\nu(c_k^{2m_k-2})\in \{2N-6, 2N-7, 2N-8\}$.

\medskip

We continue the proof by distinguishing three cases.

\medskip

{\bf Case 1:} $i(c_k^{2m_k -2})+\nu(c_k^{2m_k -2})= 2N-6 $.

In this case, by (\ref{3.3.0}) and (\ref{3.28.0}), we know that $p_{k,+}+q_{k,-}=3$ and $ i(c_k^2)=9$. It follows from  (\ref{2.18}) and (\ref{3.2.0}) that $i(c_k^2)=2i(c)\in 2\N$ since $p_{k,+}+q_{k,-}=3$, thus Case 1 cannot happen.

\medskip

{\bf Case 2:} $ i(c_k^{2m_k -2})+\nu(c_k^{2m_k -2})= 2N-7 $.

In this case, by (\ref{3.3.0}) and (\ref{3.28.0}), one of the following cases may happen.

(i)  $ i(c_k^2)=10 $ and $p_{k,+}+q_{k,-}=3$.

(ii) $ i(c_k^2)=9 $, $p_{k,+}+q_{k,-}=2 $ and $ p_{k,-}+q_{k,+}=0 $.

For (i), by  (\ref{2.18}), we have $i(c_k^2)=2i(c_k)$ and $ \hat{i}(c_k)=i(c_{k}) $ since $p_{k,+}+q_{k,-}=3$. However, there holds $ \hat i(c_{k})>5 $ by  (\ref{3.4.0}). So we have $i(c_k^2)=2i(c_k)=2\hat{i}(c_k)>10$, thus (i) of Case 2 cannot happen.

For (ii), by  (\ref{3.2.0}), there holds
\bea
p_{k,0}+q_{k,0}+r_{k,1}+r_{k,2}+h_k=1.\lb{3.29.0}
\eea
It follows from (\ref{2.18}) and  (\ref{1.2}) that
\bea
i(c_k^2)&=&2i(c_k)+p_{k,0}-q_{k,0}-3(r_{k,1}+r_{k,2})+2\sum_{j=1}^{r_{k,1}}E\left(\frac{\theta_{k,j}}{\pi}\right)+2\sum_{j=1}^{r_{k,2}}E\left(\frac{\tilde{\theta}_{k,j}}{\pi}\right)\nn\\
&>&2i(c_k)+p_{k,0}-q_{k,0}-3(r_{k,1}+r_{k,2})+2\sum_{j=1}^{r_{k,1}}\frac{\theta_{k,j}}{\pi}+2\sum_{j=1}^{r_{k,2}}\frac{\tilde{\theta}_{k,j}}{\pi}.\lb{3.30.0}
\eea
Combining (\ref{3.30.0}) and (\ref{3.29.0}), it yields
\bea
\hat i(c_{k})&=& i(c_{k})+p_{k,0}-r_{k,1}-r_{k,2}+\sum_{j=1}^{r_{k,1}}
\frac{\theta_{k,j}}{\pi}+\sum_{j=1}^{r_{k,2}}
\frac{\tilde  \theta_{k,j}}{\pi}\nn\\
&<&\frac{1}{2}\left(i(c_k^2)+p_{k,0}+q_{k,0}+r_{k,1}+r_{k,2}\right)\le \frac{1}{2}\left(i(c_k^2)+1\right),\lb{3.31.0}
\eea
which, together with  (\ref{3.4.0}), yields $i(c_k^2)>9$, thus Case (ii) cannot happen.

\medskip

{\bf Case 3:} $ i(c_k^{2m_k -2})+\nu(c_k^{2m_k -2})= 2N-8 $.

In this case, by (\ref{3.3.0}) and (\ref{3.28.0}), one of the following cases may happen.

(i)  $ i(c_k^2)=11 $ and $p_{k,+}+q_{k,-}=3$.

(ii) $ i(c_k^2)=10 $,  $p_{k,+}+q_{k,-}=2 $ and $ p_{k,-}+q_{k,+}=0 $.

(iii) $ i(c_k^2)=9 $, $p_{k,+}+q_{k,-}=2 $ and $ p_{k,-}+q_{k,+}=1 $.

(iv) $ i(c_k^2)=9 $, $p_{k,-}+q_{k,+}=0 $ and $p_{k,+}+q_{k,-}=1 $.

For (i), similar to the arguments in Case 1, it can be shown that this case cannot happen.

For (ii), similar to (\ref{3.29.0}) and the first equality in (\ref{3.30.0}), we have
\bea
&& p_{k,0}+q_{k,0}+r_{k,1}+r_{k,2}+h_k=1, \lb{3.32.0}\\
&& i(c_k^2)=p_{k,0}+q_{k,0}+r_{k,1}+r_{k,2}\quad (\mod\ 2),\lb{3.33.0}
\eea
which yields
\bea
i(c_k^2)=1+h_k\quad (\mod\ 2).\lb{3.34.0}
\eea
Then we get $h_k=1$ by $ i(c_k^2)=10$ and (\ref{3.32.0}). By  (\ref{2.18}) and (\ref{3.2.0}), we have $i(c_k^2)=2i(c_k)$ and $ \hat{i}(c_k)=i(c_{k}) $ since $p_{k,+}+q_{k,-}=2$ and $h_k=1$. However $ \hat i(c_{k})>5 $ by  (\ref{3.4.0}), thus Case (ii) cannot happen.

For (iii), by (\ref{2.18}) and (\ref{3.2.0}), we have $i(c_k^2)=2i(c_k)+p_{k,-}-q_{k,+}$ and $ \hat{i}(c_k)=i(c_k)+p_{k,-}$ since $ p_{k,-}+q_{k,+}=1 $ and $p_{k,+}+q_{k,-}=2 $. Then we obtain $2\hat{i}(c_k)=i(c_k^2)+1=10$. However $ \hat i(c_{k})>5 $ by  (\ref{3.4.0}), thus Case (iii) cannot happen.

For (iv), similar to the proof of (\ref{3.33.0}) in (ii) of Case 3, we have
\bea
i(c_k^2)=p_{k,0}+q_{k,0}+r_{k,1}+r_{k,2}\quad (\mod\ 2).\lb{3.35.0}
\eea
Therefore by (\ref{3.2.0}), (\ref{3.35.0}), $i(c_k^2)=9$ and $p_{k,+}+q_{k,-}=1$, we get
\bea
p_{k,0}+q_{k,0}+r_{k,1}+r_{k,2}=1.\lb{3.36.0}
\eea
Then similar to the proof of (\ref{3.31.0}) in (ii) of Case 2,  we have
\bea
\hat i(c_{k})<\frac{1}{2}\left(i(c_k^2)+p_{k,0}+q_{k,0}+r_{k,1}+r_{k,2}\right)= \frac{1}{2}\left(i(c_k^2)+1\right),\lb{3.37.0}
\eea
which, together with  (\ref{3.4.0}), yields $i(c_k^2)>9$, thus Case (iv) cannot happen.

This completes the proof of Lemma 3.2.\hfill\hb

\medskip

Under the assumption (FCG), Theorem 1.1 of \cite{Dua2} shows that there exist at least two elliptic closed geodesics $c_1$ and $c_2$ on $(S^4,F)$ whose flag curvature satisfies $\left(\frac{\lm}{1+\lm}\right)^2<K\le 1$. The following Lemma gives some properties of these two closed geodesics which will be useful in the proof of Theorem 1.1.

\medskip

{\bf Lemma 3.3.} (cf. Lemma 3.1 and Lemma 3.3 of \cite{Dua1} and Section 3 of \cite{Dua2}) {\it Under the assumption (FCG), there exist at least two prime elliptic closed geodesics $c_1$ and $c_2$ on $(S^4,F)$ whose flag curvature satisfies $\left(\frac{\lm}{1+\lm}\right)^2<K\le 1$. Moreover, there exist infinitely many pairs of $(q+1)$-tuples $(N, m_1, m_2,\cdots, m_q)\in\N^{q+1}$ and $(N', m_1', m_2',\cdots, m_q')\in\N^{q+1}$ such that
\bea && i(c_1^{2m_1})+\nu(c_1^{2m_1})=2N+3,\quad \ol{C}_{2N+3}(E,c_1^{2m_1})=\Q,\lb{3.38.0}\\
&& i(c_2^{2m_2'})+\nu(c_2^{2m_2'})=2N'+3,\quad \ol{C}_{2N'+3}(E,c_2^{2m_2'})=\Q,\lb{3.39.0}\\
&& p_{k,-}=q_{k,+}=r_{k,3}=r_{k,4}=r_{k,6}=h_k=0,\quad k=1,2,\lb{3.40.0}\\
&& r_{1,2}=\Delta_1\ge 1,\qquad r_{2,2}=\Delta_2'\ge 1, \lb{3.41.0}\\
&& \Delta_k+\Delta_k'=r_{k,2},\quad k=1,2,\lb{3.42.0}\eea
where we can require $3|N$ or $3|N'$ as remarked in Theorem 2.8 and \bea \Delta_k' \equiv \sum_{0<\{m_k'\th/\pi\}<\delta}S^-_{M_k}(e^{\sqrt{-1}\th}),\quad k=1,2.\lb{3.43.0}\eea

In addition, for these two closed geodesics $c_1$ and $c_2$, there holds
\bea k_{\nu(c_k^{n(c_k)})}^{\ep(c_k^{n(c_k)})}(c_k^{n(c_k)})=1,\quad k_{j}^{\ep(c_k^{n(c_k)})}(c_k^{n(c_k)})=0,\quad\forall\  0\le j<\nu(c_k^{n(c_k)}),\ k=1,2. \lb{3.44.0}\eea}

\subsection{Proof of Theorem 1.1}

In this section, let $(S^4,F)$ be a Finsler sphere of dimension $4$ with its reversibility $\lambda$ and flag curvature $K$ satisfying $\frac{25}{9}\left(\frac{\lm}{1+\lm}\right)^2<K\le 1$ with $\lambda<\frac{3}{2}$. In order to prove Theorem 1.1, according to Lemma 3.3 and Theorem 1.1 of \cite{Dua1}, we make the following assumption.

\smallskip

{\bf (TCG)} {\it Suppose that there exist exactly two prime elliptic closed geodesics $c_1$ and $c_2$ possessing all properties listed
in Lemma 3.3, and the third prime closed geodesic $c_3$ on such $(S^4,F)$.}

\smallskip

In order to count the contribution of $c_k^m$ to the Morse-type number $M_q$, for the sake of convenience, we set
\bea
M_q(k,m)=\dim{\ol{C}}_q(E, c^m_k),\quad \forall\ 1\le k\le 3,\ m\ge1,\ q\in\N_0. \lb{4.0.1}
\eea

Next we fix $\bar{m}=4$. Before proving Theorem 1.1, firstly we establish several crucial lemmas.

\medskip

{\bf Lemma 3.4.} {\it For an integer $q$ satisfying $2N-9\le q\le 2N+17$, there holds}
\bea
M_q=\left\{\matrix{
    \sum_{1\le k\le 3, \atop 1\le m\le 2}M_q(k,2m_k-m),\quad &&{\it if}\quad q=2N-9,  \cr
    \sum_{1\le k\le 3}M_q(k,2m_k-1),\quad &&{\it if}\quad 2N-8\le q\le 2N-4, \cr
    \sum_{1\le k\le 3, \atop 0\le m\le 1}M_q(k,2m_k-m),\quad &&{\it if}\quad q=2N-3,  \cr
    \sum_{1\le k\le 3}M_q(k,2m_k),\quad &&{\it if}\quad 2N-2\le q\le 2N+2,  \cr
    \sum_{1\le k\le 3, \atop 0\le m\le 1}M_q(k,2m_k+m),\quad &&{\it if}\quad q=2N+3,  \cr
    \sum_{1\le k\le 3}M_q(k,2m_k+1),\quad &&{\it if}\quad 2N+4\le q\le 2N+8, \cr
    \sum_{1\le k\le 3, \atop 1\le m\le 2}M_q(k,2m_k+m),\quad &&{\it if}\quad  2N+9\le q\le 2N+14,  \cr
    \sum_{1\le k\le 3, \atop 1\le m\le 3}M_q(k,2m_k+m),\quad &&{\it if}\quad  2N+15\le q\le 2N+17.  \cr}\right.
\eea

\medskip

{\bf Proof.} According to Lemma 2.1, (\ref{2.4}) and (i) of Lemma 2.2, we have
\bea
M_q&=&\sum_{1\le k\le 3, \atop m\ge 1}M_q(k,m)=\sum_{1\le k\le 3, \atop m\ge 1}k_{q-i(c_k^m)}^{\ep(c_k^m)}(c_k^m)\nn\\
&=&\sum_{1\le k \le 3, \atop m\in\left\{m\in\N|i(c_k^m)\le q \le i(c_k^m)+\nu(c_k^m)\right\}}k_{q-i(c_k^m)}^{\ep(c_k^m)}(c_k^m)=\sum_{1\le k \le 3, \atop m\in\left\{m\in\N |i(c_k^m)\le q \le i(c_k^m)+\nu(c_k^m)\right\}}M_q(k,m).\lb{3.47.0}
\eea

On one hand, by (\ref{3.19.0}), (\ref{3.15.0}), (\ref{3.24.0}), (\ref{3.23.0}) and (\ref{3.17.0}), it yields
\bea
i(c_k^m)+\nu(c_k^m)\le\left\{\matrix{
    2N-15,&&\quad {\it if}\quad 1\le m\le 2m_k-4,  \cr
    2N-12,&&\quad {\it if}\quad m=2m_k-3,  \cr
    2N-9,&&\quad {\it if}\quad m=2m_k-2,  \cr
    2N-3,&&\quad {\it if}\quad m=2m_k-1,  \cr
    2N+3,&&\quad {\it if}\quad m=2m_k.  \cr}\right.\lb{3.48.0}
\eea

On the other hand, by (\ref{3.16.0}), (\ref{3.18.0}) and (\ref{3.20.0}), it yields
\bea
i(c_k^m)\ge\left\{\matrix{
    2N-3,&&\quad {\it if}\quad m=2m_k,  \cr
    2N+3,&&\quad {\it if}\quad m=2m_k+1,  \cr
    2N+9,&&\quad {\it if}\quad m=2m_k+2,  \cr
    2N+15,&&\quad {\it if}\quad m=2m_k+3,  \cr
    2N+18,&&\quad {\it if}\quad m\ge 2m_k+4.  \cr}\right.\lb{3.49.0}
\eea
Combining (\ref{3.47.0})-(\ref{3.49.0}), we get Lemma 3.4. \hfill\hb

\medskip

{\bf Lemma 3.5.} {\it For some tuple $(k,m)$ with $k=1,2,3$ and $m\in\N$, if there exist some integers $q_1, q_2\in\N$ satisfying
\bea
q_1\le i(c_k^m) \quad \mbox{and} \quad i(c_k^m)+\nu(c_k^m)\le q_2,\nn
\eea
then there holds
\bea
M_{q_1}(k,m)+M_{q_2}(k,m)\le 1.
\eea
Furthermore, if  $M_{q_1}(k,m)+M_{q_2}(k,m)=1$, then
\bea
M_q(k,m)=0,\quad \forall\  q\neq q_1,q_2.
\eea}

{\bf Proof.} This follows directly from Lemma 2.1, (\ref{2.4}), (i) and (ii) of Lemma 2.2. \hfill\hb

\medskip

{\bf Lemma 3.6.} {\it For some $k\in \{1,2,3\}$, assume that either $M_{2N-q\pm1}(k,2m_k-1)\ge 1$
or $M_{2N+q\pm1}(k,2m_k+1)\ge 1$ for some even $q\in \N$, then there
exists a continuous path $f_k\in C([0,1],\Om^0(P_{c_k}))$ such that $f_k(0)=P_{c_k}$ and $f_k(1)$ belongs to one of the following five cases:

(i) $I_4\,\dm\,H(2)$,

(ii) $N_1(1,1)\,\dm\,I_2\,\dm\,N_1(1,-1)$,

(iii) $I_4\,\dm\,N_1(1,-1)$,

(iv) $N_1(1,1)\,\dm\, I_4$,

(v) $I_6$.

And in either case, the index iteration formula of $c_k^m$ can be written as follows:
\bea
i(c_k^m)=qm-p_{k,-}-p_{k,0},\quad \forall\ m\ge 1. \lb{3.52.0}
\eea }

{\bf Proof.} We only give the proof under the assumption $M_{2N-q\pm1}(k,2m_k-1)\ge 1$. The proof under the assumption $M_{2N+q\pm1}(k,2m_k+1)\ge 1$
is similar.

First, by Lemma 3.5 and the assumption $M_{2N-q\pm1}(k,2m_k-1)\ge 1$ in Lemma 3.6, we have
\bea
i(c_k^{2m_k-1}) \le 2N-q-2, \qquad i(c_k^{2m_k-1})+\nu(c_k^{2m_k-1}) \ge 2N-q+2, \lb{3.53.0}
\eea
which, together with $\nu(c_k^{2m_k-1})=\nu(c_k)$ by (\ref{2.22}), implies $\nu(c_k)=p_{k,-}+2p_{k,0}+p_{k,+}\in\{4,5,6\}$.

If $\nu(c_k)=4$, by (\ref{3.53.0}), we have $i(c_k^{2m_k-1})= 2N-q-2$ and $i(c_k^{2m_k-1})+\nu(c_k^{2m_k-1})=2N-q+2$. And we also have $p_{k,-}+p_{k,0}+p_{k,+}=2$ or $3$ by (\ref{3.2.0}). Since $i(c_k^{2m_k-1})= 2N-q-2\in 2\N$, we get $i(c_k)\in 2\N$ by (\ref{2.24}), then by Proposition 2.7 and the symplectic additivity of symplectic paths (cf. Theorem 6.2.6 of \cite{Lon3}), we must have $p_{k,+}+h_{k,+}=1$. Therefore, if $p_{k,-}+p_{k,0}+p_{k,+}=2$, we must have $p_{k,0}=2$ since $\nu(c_k)=4$ and $h_{k,+}=1$, i.e., $f_k(1)=I_4\,\dm\,H(2)$.  If $p_{k,-}+p_{k,0}+p_{k,+}=3$, we must have $p_{k,0}=1$, $p_{k,+}=1$ and $p_{k,-}=1$, i.e., $f_k(1)=N_1(1,1)\,\dm\,I_2\,\dm\,N_1(1,-1)$.

If $\nu(c_k)=5$, we have $p_{k,-}+p_{k,0}+p_{k,+}=3$. Then we must have $p_{k,0}=2$
and $p_{k,-}+p_{k,+}=1$. So either $i(c_k^{2m_k-1})= 2N-q-2$ when $p_{k,+}=1$, or $i(c_k^{2m_k-1})= 2N-q-3$ when $p_{k,-}=1$ by Proposition 2.7, the symplectic additivity and (\ref{3.53.0}), i.e., $f_k(1)=I_4\,\dm\,N_1(1,-1)$ or $N_1(1,1)\,\dm\, I_4$.

If $\nu(c_k)=6$, we must have $p_{k,0}=3$, and then we have $i(c_k^{2m_k-1})= 2N-q-3$ by Proposition 2.7, the symplectic additivity and (\ref{3.53.0}), i.e., $f_k(1)=I_6$.

Note that by (\ref{2.24}), (\ref{3.10.0}) and (\ref{3.12.0}), we get
$$
i(c_k^{2m_k-1})=2N-i(c_k)-2p_{k,-}-2p_{k,0},
$$
and by the above arguments in either case, we have $i(c_k)=q-p_{k,-}-p_{k,0}$.  Then by (\ref{2.18}), we have
$$
i(c_k^m)=qm-p_{k,-}-p_{k,0}.
$$
This completes the proof of Lemma 3.6. \hfill\hb

\medskip

{\bf Proof of Theorem 1.1.}

\medskip

At first, we consider the contribution of $c_k^m$ with $k=1,2,3$ and $m\in\N$ to the Morse-type numbers in Claim 1 and Claim 2. And then we use Claim 3 and Claim 4 to complete the proof of Theorem 1.1.

\medskip

{\bf Claim 1:} {\it For $2N-2\le q\le 2N+2$, there holds: (i) $M_q(1,m)=0$ for any $m\in\N$, and (ii) $M_q(k,m)=0$ for $k=2,3$ and $m\neq 2m_k$. In addition,
$M_q(2,2m_2)=0$ for $q=2N-2,2N,2N+2$ and $M_{2N-1}(2,2m_2)+M_{2N+1}(2,2m_2)\le 1$.}

 \medskip

{\bf Proof.} By Lemma 3.4, we know that $M_q(k,m)=0$
for $2N-2\le q\le 2N+2$, $m \neq 2m_k$ and $k=1,2,3$.
 Also note that by (\ref{3.38.0}),
(\ref{3.44.0}), Lemma 2.1 and (\ref{2.4}), we have
$$
M_q(1,2m_1)=0\quad \mbox{for}\quad 2N-2\le q\le 2N+2.
$$

On one hand, there holds
$$\nu(c_2^{2m_2})=\nu(c_2^{2m'_2})$$
by the choices of $m_2$ and $m_2'$ in (\ref{2.27}) of Theorem 2.8. On the other hand, it yields
$$i(c_2^{2m_2} ) = i(c_2^{2m_2'} ) \quad(\mod\ 2) $$
by (\ref{2.18}) of Theorem 2.6. So, $i(c_2^{2m_2} ) + \nu(c_2^{2m_2} )$ is odd since $i(c_2^{2m'_2} ) + \nu(c_2^{2m'_2} )$ is odd by (\ref{3.39.0}) of Lemma 3.3, which implies that $M_q(2,2m_2)=0$ for $q=2N-2, 2N, 2N+2$
and $M_{2N-1}(2,2m_2)+M_{2N+1}(2,2m_2)\le 1$ by (\ref{3.44.0}), Lemma 2.1 and (\ref{2.4}). Hence, Claim 1 holds.

\medskip

{\bf Claim 2:} {\it $M_{2N-1}(2,2m_2)=M_{2N+1}(2,2m_2)=0$.}

\medskip

{\bf Proof.} Otherwise, by Claim 1 we have
\bea
M_{2N-1}(2,2m_2)+M_{2N+1}(2,2m_2)=1. \lb{3.54.0}
\eea
Then by (\ref{3.16.0}) and Lemma 3.5, it yields
\bea
M_{2N-3}(2,2m_2)=0. \lb{3.55.0}
\eea
By (\ref{2.8}) and Theorem 2.5, we have $M_{2N-1}\ge b_{2N-1}=1$ and $M_{2N+1}\ge b_{2N+1}=1$, then we get $M_{2N-1}(3,2m_3)+M_{2N+1}(3,2m_3)\ge 1$  by Lemma 3.4, (i) of Claim 1 and (\ref{3.54.0}). Thus by (\ref{3.16.0}), (\ref{3.17.0}) and Lemma 3.5, it yields
\bea
M_{2N-3}(3,2m_3)=M_{2N+3}(3,2m_3)=0. \lb{3.56.0}
\eea
So, by Lemma 3.4, Claim 1, (\ref{3.54.0}) and (\ref{3.56.0}), we obtain
\bea
\sum_{q=2N-2}^{2N+2}(-1)^q M_q=\sum_{q=2N-2}^{2N+2}(-1)^q \sum_{1\le k \le 3}M_q(k,2m_k)=\sum_{q=2N-3}^{2N+3}(-1)^q M_q(3,2m_3)-1.\lb{3.57.0}
\eea
And by (\ref{3.16.0}), (\ref{3.17.0}), Lemma 2.1 and (\ref{2.4}), we have
\bea
\sum_{q=2N-3}^{2N+3}(-1)^q M_q(3,2m_3)
=\sum_{0\le l \le 6}(-1)^{i(c_3^{2m_3})+l}k_{l}^{\ep(c_3^{2m_3})}(c_3^{2m_3}).\lb{3.58.0}
\eea

On the other hand, by (\ref{2.8}) and Theorem 2.5, we have
\bea
\sum_{q=2N-2}^{2N+2}(-1)^q M_q\ge \sum_{q=2N-2}^{2N+2}(-1)^q b_q=-2.\lb{3.59.0}
\eea
Combining (\ref{3.57.0})-(\ref{3.59.0}), we get
\bea
\chi(c_3^{2m_3})=\sum_{0\le l \le 6}(-1)^{i(c_3^{2m_3})+l}k_{l}^{\ep(c_3^{2m_3})}(c_3^{2m_3}) \ge -1.\lb{3.60.0}
\eea
Note that since $n(c_3)|2m_3$ and $\nu(c_3^{2m_3}) = \nu(c_3^{n(c_3)})$ by (\ref{2.7}) and (\ref{2.27}), there holds
\bea
\chi(c_3^{n(c_3)})=\chi(c_3^{2m_3})\ge -1.\lb{3.61.0}
\eea

By (\ref{3.38.0}), (\ref{3.44.0}), Lemma 2.1 and (\ref{2.4}), it yields $M_{2N-3}(1,2m_1)=0$. Together with (\ref{3.55.0}) and (\ref{3.56.0}), we have
\bea
\sum_{1\le k\le 3}M_{2N-3}(k,2m_k)=0.\lb{3.62.0}
\eea
Then combining Lemma 3.4 and (\ref{3.62.0}), we obtain that
\bea
M_{2N-3}=\sum_{1\le k\le 3}M_{2N-3}(k,2m_k-1).\lb{3.63.0}
\eea
One one hand, by (\ref{3.23.0}) and Lemma 3.5 it yields
\bea
M_{2N-3}(k,2m_k-1)\le 1,\quad \forall\ k=1,2,3,\lb{3.64.0}
\eea
then it follows from (\ref{3.63.0}) and (\ref{3.64.0}) that  $M_{2N-3}\le 3$. On the other hand, we have  $M_{2N-3}\ge b_{2N-3}=2$ by (\ref{2.8}) and Theorem 2.5. So it yields $M_{2N-3}\in \{2,3\}$.

We continue the proof by distinguishing two cases.

\medskip

{\bf Case 1:} $M_{2N-3}=3$.

\medskip

In this case, it follows from (\ref{3.63.0}) and (\ref{3.64.0}) that   $M_{2N-3}(k,2m_k-1)=1$ for $k=1,2,3$. Then according to (\ref{3.23.0}) and Lemma 3.5, there holds
\bea
M_{2N-5}(k,2m_k-1)=0, \quad\forall\  k=1,2,3. \lb{3.65.0}
\eea
Combining (\ref{3.65.0}) and Lemma 3.4, we get $M_{2N-5}=0$. But by (\ref{2.8}) and Theorem 2.5, we have $M_{2N-5} \ge b_{2N-5}=1$. This is a contradiction.

\medskip

{\bf Case 2:} $M_{2N-3}=2$.

\medskip

In this case, it follows from (\ref{3.63.0}) and (\ref{3.64.0}) that one of the following cases may happen:

(i) $M_{2N-3}(3,2m_3-1)=0$ and $M_{2N-3}(1,2m_1-1)=M_{2N-3}(2,2m_1-1)=1$.

(ii) $M_{2N-3}(3,2m_3-1)=1$ and $M_{2N-3}(1,2m_1-1)+M_{2N-3}(2,2m_2-1)=1$.

For (i), by (\ref{3.23.0}) and Lemma 3.5, there holds
\bea
M_q(k,2m_k-1)=0,\quad \forall\ q\neq 2N-3,\ k=1,2.\lb{3.67.0}
\eea
So, according to Lemma 3.4 and (\ref{3.67.0}), we have
\bea
M_{2N-5}&=&\sum_{1\le k\le 3}M_{2N-5}(k,2m_k-1)=M_{2N-5}(3,2m_3-1),\lb{3.68.0}\\
M_{2N-7}&=&\sum_{1\le k\le 3}M_{2N-7}(k,2m_k-1)=M_{2N-7}(3,2m_3-1).\lb{3.69.0}
\eea
By (\ref{2.8}) and Theorem 2.5, we have $M_{2N-5} \ge b_{2N-5}=1$ and $M_{2N-7} \ge b_{2N-7}=1$, then it follows from (\ref{3.68.0}) and (\ref{3.69.0}) that
\bea
M_{2N-5}(3,2m_3-1)\ge 1,\quad M_{2N-7}(3,2m_3-1)\ge 1. \lb{3.70.0}
\eea
So the assumption with $q=6$ in Lemma 3.6 is satisfied, and then by (\ref{3.52.0}) and (\ref{2.7}), we have
\bea
i(c_3^m) &=& 6m-p_{3,-}-p_{3,0},\lb{3.71.0}\\
n(c_3) &=& 1.\lb{3.72.0}
\eea

Notice that $\nu(c_3^{2m_3-1})=\nu(c_3^{2m_3-2})=\nu(c_3)$ by (\ref{2.19}) for each of five cases in Lemma 3.6. Together with Lemma 2.1, (\ref{2.4}), (iv) of Lemma 2.2, (\ref{3.71.0}), (\ref{3.72.0}) and (i) of Case 2,  it yields
\bea
M_{2N-9}(3,2m_3-2)&=&k_{2N-9-i(c_3^{2m_3-2})}^{\ep(c_3^{2m_3-2})}(c_3^{2m_3-2})\nn\\
&=&k_{2N-9-i(c_3^{2m_3-2})}^{\ep(c_3^{2m_3-1})}(c_3^{2m_3-1})\nn\\
&=&M_{2N-9-i(c_3^{2m_3-2})+i(c_3^{2m_3-1})}(3,2m_3-1)\nn\\
&=&M_{2N-3}(3,2m_3-1)=0.\lb{3.73.0}
\eea

Through comparing (\ref{3.9.0}) and (\ref{3.71.0}), we get $2N=12m_3$. Then by (\ref{3.71.0}) and (\ref{3.2.0}) it yields
\bea
i(c_3^{2m_3-1})=12m_3-6-p_{3,-}-p_{3,0}\ge 2N-9,\lb{3.74.0}
\eea
Then by (\ref{3.70.0}) and Lemma 3.5, there holds
\bea
M_{2N-9}(3,2m_3-1)=0.\lb{3.75.0}
\eea
Finally, combining (\ref{3.67.0}), (\ref{3.73.0}), (\ref{3.75.0}) and Lemma 3.4, we obtain
\bea
M_{2N-9}=\sum_{1\le k\le 3, \atop 1\le m\le 2}M_{2N-9}(k,2m_k-m)=M_{2N-9}(1,2m_1-2)+M_{2N-9}(2,2m_2-2).\lb{3.76.0}
\eea

By (\ref{3.24.0}) and (ii) of Lemma 2.2, we get  $M_{2N-9}(k,2m_k-2)\le 1$ for $k=1,2$. Then it follows from (\ref{3.76.0}) that $M_{2N-9}\le 2$. However by (\ref{2.8}) and Theorem 2.5, we have $M_{2N-9}\ge b_{2N-9}=2$. Thus we get $M_{2N-9}= 2$.

Now, in this case, we have
\bea
M_{2N-3}=b_{2N-3},\quad  M_{2N-9}=b_{2N-9}.\lb{3.77.0}
\eea
Combining (\ref{2.8}), Theorem 2.5 and (\ref{3.77.0}), we obtain
\bea
\sum_{q=2N-8}^{2N-4} (-1)^q M_{q}=\sum_{q=2N-8}^{2N-4} (-1)^q b_{q}=-2.\lb{3.78.0}
\eea
Note that for $2N-8\le q\le 2N-4$, it follows from (\ref{3.23.0}) and (i) of this case that $M_q(k,2m_k-1)=0$ for $k=1,2$. So we get $M_q=M_q(3,2m_3-1)$ by Lemma 3.4, and then it yields
\bea
\sum_{q=2N-8}^{2N-4} (-1)^q M_{q}=\sum_{q=2N-8}^{2N-4} (-1)^q M_{q}(3,2m_3-1).\lb{3.79.0}
\eea
Note that  $2N-9\le i(c_3^{2m_3-1})\le i(c_3^{2m_3-1})+\nu(c_3^{2m_3-1}) \le 2N-3$ by (\ref{3.74.0}) and (\ref{3.23.0}), according to (\ref{3.75.0}), (i) of this case, Lemma 2.1 and (\ref{2.4}), we obtain
\bea
\sum_{q=2N-8}^{2N-4} (-1)^q M_{q}(3,2m_3-1)&=&\sum_{q=2N-9}^{2N-3} (-1)^q M_{q}(3,2m_3-1)\nn\\
&=&\sum_{0\le l \le 6}(-1)^{i(c_3^{2m_3-1})+l}k_{l}^{\ep(c_3^{2m_3-1})}(c_3^{2m_3-1}).\lb{3.80.0}
\eea
Combining (\ref{3.78.0})-(\ref{3.80.0}), we get
\bea
\chi(c_3^{2m_3-1})=\sum_{0\le l \le 6}(-1)^{i(c_3^{2m_3-1})+l}k_{l}^{\ep(c_3^{2m_3-1})}(c_3^{2m_3-1})=-2.\lb{3.81.0}
\eea
However, since $n(c_3)=1$ by (\ref{3.72.0}), it follows from (iv) of Lemma 2.2 and (\ref{3.81.0}) that $\chi(c_3)=\chi(c_3^{2m_3-1})=-2$, which contradicts to (\ref{3.61.0}), thus Case (i) cannot happen.

For (ii), without loss of generality, we assume that $M_{2N-3}(1,2m_1-1)=0$. Then by (\ref{3.23.0}) and  Lemma 3.5, $M_q(k,2m_k-1)=0$ for $q\neq 2N-3$ and $k=2,3$. So, according to Lemma 3.4,
we have
\bea
M_{2N-5}=M_{2N-5}(1,2m_1-1),\quad M_{2N-7}=M_{2N-7}(1,2m_1-1). \lb{3.82.0}
\eea
By (\ref{2.8}) and Theorem 2.5, we have $M_{2N-5} \ge b_{2N-5}=1$ and $M_{2N-7} \ge b_{2N-7}=1$, then it follows from (\ref{3.82.0}) that $M_{2N-5}(1,2m_1-1)\ge 1$ and $M_{2N-7}(1,2m_1-1)\ge 1$. Thus by Lemma 3.6, there holds $f_1(0)=P_{c_1}$ and $f_1(1)$ belongs to one of five cases in Lemma 3.6, which contradicts (\ref{3.40.0}) in Lemma 3.3.

This completes the proof of Claim 2.

\medskip

{\bf Claim 3:} {\it $c_1$ and $c_2$ are irrationally elliptic.}

\medskip

{\bf Proof.} By (\ref{3.41.0}) and (\ref{3.42.0}), there holds $\Delta_2 = 0$. Then, together with the fact that $r_{2,3}=r_{2,4} = 0$ from (\ref{3.40.0}), it follows from (\ref{3.41.0}) and (\ref{3.21.0}) that
\bea
2N+1 &\ge& i(c_2^{2m_2}) + \nu(c_2^{2m_2} )\lb{3.84.0} \\
&=&2N+(p_{2,0} +p_{2,+} +q_{2,-} +q_{2,0} +r_{2,1}+2r_{2,5}-r_{2,2})\nn\\
&\ge& 2N - 3, \lb{3.85.0}
\eea
where (\ref{3.84.0}) holds by the fact that $p_{2,0} +p_{2,+} +q_{2,-} +q_{2,0} +r_{2,1}+2r_{2,5}\le 2$ from (\ref{3.2.0}) and (\ref{3.41.0}), and $r_{2,2}\ge 1$ from (\ref{3.41.0}), and the equality in (\ref{3.85.0}) holds if and only if  $r_{2,2}=3$.

On the other hand, by Claim 2, we have $i(c_2^{2m_2} ) + \nu(c_2^{2m_2} ) \notin \{2N-1,  2N+1\}$. Note that $i(c_2^{2m_2} )+\nu(c_2^{2m_2} )=i(c_2^{2m'_2} )+\nu(c_2^{2m'_2})=1\ (\mod\ 2)$ by (\ref{3.39.0}). Thus, by (\ref{3.84.0}), we obtain $i(c_2^{2m_2} ) + \nu(c_2^{2m_2} ) \le 2N-3$, which together with (\ref{3.85.0}) implies $r_{2,2} = 3$, i.e., $c_2$ is irrationally elliptic. By the symmetric properties of $c_1$ and $c_2$ in Lemma 3.3 (or, more precisely, replacing $N$ with $N'$ in the above arguments), we conclude that $c_1$ is also irrationally elliptic.

This completes the proof of Claim 3.

\medskip

{\bf Claim 4:} {\it $c_3$ is non-hyperbolic.}

\medskip

{\bf Proof.} Assume that $c_3$ is hyperbolic, which, together with the assumption (TCG) and Claim 3, implies the Finsler metric $F$ on $S^4$ is bumpy. Then by Theorem 1.1 in \cite{DLW}, we know that there exist at least four distinct non-hyperbolic prime closed geodesics, which contradicts the assumption (TCG). Thus  by $c_3$ is non-hyperbolic.

Therefore, by the assumption (TCG), Claims 3 and 4 complete the proof of Theorem 1.1.   \hfill\hb

\setcounter{figure}{0}
\setcounter{equation}{0}
\section{Some further information about the third closed geodesic}  

In this section, under the assumption (TCG), we further study the third closed geodesic $c_3$ and obtain some much precise information about it (cf. Theorem 4.2 below). At first, we establish a similar result as Lemma 3.6.

\medskip

{\bf Lemma 4.1.} {\it For some $k\in \{1,2,3\}$, assume that either $M_{2N-q\pm1}(k,2m_k-2)\ge 1$
or $M_{2N+q\pm1}(k,2m_k+2)\ge 1$ for some even $q\in \N$, then there
exists a continuous path $f_k\in C([0,1],\Om^0(P_{c_k}))$ such that $f_k(0)=P_{c_k}$ and $f_k(1)$ belongs to one of the following cases:

(i) $I_{2p_{k,0}}\,\dm\,(-I_{2q_{k,0}})\,\dm\,H(2)$ with $p_{k,0}+q_{k,0}=2$,

(ii) $N_1(1,1)^{\dm p_{k,-}}\,\dm\,N_1(-1,-1)^{\dm q_{k,+}}\,\dm\,I_{2p_{k,0}}\,\dm\,(-I_{2q_{k,0}})\,\dm\,N_1(1,-1)^{\dm p_{k,+}}\,\dm\,N_1(-1,1)^{\dm q_{k,-}}$ with $p_{k,-}+q_{k,+}=1$, $p_{k,0}+q_{k,0}=1$ and $p_{k,+}+q_{k,-}=1$,

(iii) $I_{2p_{k,0}}\,\dm\,(-I_{2q_{k,0}})\,\dm\,N_1(1,-1)^{\dm p_{k,+}}\,\dm\,N_1(-1,1)^{\dm q_{k,-}}$ with $p_{k,0}+q_{k,0}=2$ and $p_{k,+}+q_{k,-}=1$,

(iv) $N_1(1,1)^{\dm p_{k,-}}\,\dm\,N_1(-1,-1)^{\dm q_{k,+}}\,\dm\,I_{2p_{k,0}}\,\dm\,(-I_{2q_{k,0}})$ with $p_{k,0}+q_{k,0}=2$ and $p_{k,-}+q_{k,+}=1$,

(v) $I_{2p_{k,0}}\,\dm\,(-I_{2q_{k,0}})$ with $p_{k,0}+q_{k,0}=3$.

And in either case, the index iteration formula of $c_k^m$ can be written as follows:
\bea
i(c_k^m)=\frac{q}{2}m-(p_{k,-}+p_{k,0})-\frac{1+(-1)^m}{2}(q_{k,0}+q_{k,+}).
\eea}

{\bf Proof.} We only give the proof under the assumption $M_{2N+q\pm1}(k,2m_k+2)\ge 1$. The proof under the assumption $M_{2N-q\pm1}(k,2m_k-2)\ge 1$
is similar.

First, by Lemma 3.5 and the assumption $M_{2N+q\pm1}(k,2m_k+2)\ge 1$ in Lemma 4.1, we have
\bea
i(c_k^{2m_k+2})\le 2N+q-2,\quad i(c_k^{2m_k+2})+\nu(c_k^{2m_k+2})\ge 2N+q+2,\lb{5.2}
\eea
which, together with $\nu(c_k^{2m_k+2})=\nu(c_k^2)$ by (\ref{2.22}), implies $\nu(c_k^2)=p_{k,-}+2p_{k,0}+p_{k,+}+q_{k,-}+2q_{k,0}+q_{k,+}\in\{4,5,6\}$.

If $\nu(c_k^2)=4$, by (\ref{5.2}), we have $i(c_k^{2m_k+2})= 2N+q-2$ and $i(c_k^{2m_k+2})+\nu(c_k^{2m_k+2})=2N+q+2$. And by (\ref{3.2.0}) we also have
\bea
p_{k,-}+p_{k,0}+p_{k,+}+q_{k,-}+q_{k,0}+q_{k,+}\in\{2,3\}. \lb{5.3}
\eea
Since $i(c_k^{2m_k+2})= 2N+q-2\in 2\N$, we get $i(c_k^2)\in 2\N$ by (\ref{3.9.0}).
Note that by (\ref{2.18}) and (\ref{3.2.0}), we have
\bea
i(c_k^2)&=&p_{k,-}+p_{k,0}+q_{k,0}+q_{k,+}+r_{k,1}+r_{k,2}\quad(\mod\ 2)\nn\\
&=&1+p_{k,+}+q_{k,-}+h_k\quad(\mod\ 2). \lb{5.4}
\eea
So, by (\ref{5.3}) and (\ref{3.2.0}), we get $p_{k,+}+q_{k,-}+h_{k}=1$ since $\nu(c_k^2)=4$. Therefore, if $p_{k,-}+p_{k,0}+p_{k,+}+q_{k,-}+q_{k,0}+q_{k,+}=2$, we must have $p_{k,0}+q_{k,0}=2$ and $h_k=1$, if $p_{k,-}+p_{k,0}+p_{k,+}+q_{k,-}+q_{k,0}+q_{k,+}=3$, we must have $p_{k,0}+q_{k,0}=1$, $p_{k,+}+q_{k,-}=1$ and $p_{k,-}+q_{k,+}=1$.

If $\nu(c_k^2)=5$, we have $p_{k,-}+p_{k,0}+p_{k,+}+q_{k,-}+q_{k,0}+q_{k,+}=3$. Then we must have $p_{k,0}+q_{k,0}=2$
and $p_{k,-}+q_{k,+}+p_{k,+}+q_{k,-}=1$. And by (\ref{3.9.0}) and (\ref{5.4}) we have $i(c_k^{2m_k+2})= 2N+q-2$ when $p_{k,+}+q_{k,-}=1$, $i(c_k^{2m_k+2})= 2N+q-3$ when $p_{k,-}+q_{k,+}=1$.

If $\nu(c_k)=6$, we must have $p_{k,0}+q_{k,0}=3$, and then we have $i(c_k^{2m_k+2})= 2N+q-3$ by (\ref{5.4}) and (\ref{3.9.0}).

Note that by (\ref{3.9.0}), we get $i(c_k^{2m_k+2})=2N+i(c_k^2)$,
and in either case, we have
\bea
i(c_k^2)=q-p_{k,-}-p_{k,0}-q_{k,+}-q_{k,0}.\lb{5.5}
\eea
Then by (\ref{2.18}) and according to the precise cases in Lemma 4.1, we have
\bea
i(c_k^2)=2i(c_k)+p_{k,-}+p_{k,0}-q_{k,0}-q_{k,+}.\lb{5.6}
\eea
Combining (\ref{2.18}), (\ref{5.5}) and (\ref{5.6}), it yields
$$
i(c_k^m)=\frac{q}{2}m-(p_{k,-}+p_{k,0})-\frac{1+(-1)^m}{2}(q_{k,0}+q_{k,+}).
$$

This completes the proof of Lemma 4.1. \hfill\hb

\medskip

{\bf Theorem 4.2.} {\it For every Finsler metric F on $S^4$ with reversibility $\lm$ and flag curvature $K$ satisfying $\frac{25}{9}(\frac{\lm}{\lm+1})^2 < K \le 1$ with $\lambda<\frac{3}{2}$, suppose that there exist precisely three prime closed geodesics $c_1$, $c_2$ and $c_3$, then both $c_1$ and $c_2$ are irrationally elliptic with $i(c_1)=3$ and $ i(c_2)=9$, and $c_3$ is non-hyperbolic and must belong to one of the following precise classes:

(i) $i(c_3)=3$ and $P_{c_3}\approx N_1(1,1)\,\dm\, I_4$,

(ii) $i(c_3)=3$ and $P_{c_3}\approx I_6$,

(iii) $i(c_3)=4$ and $P_{c_3}\approx N_1(1,1)\,\dm\,I_2\,\dm\,N_1(1,-1)$,

(iv) $i(c_3)=4$ and $P_{c_3}\approx I_4\,\dm\,N_1(1,-1)$,

(v) $i(c_3)=4$ and $P_{c_3}\approx I_4\,\dm\,H(2)$,

\noindent where and below, ``$P_{c_k}\approx A$'' means that there exists a continuous path $f_k\in C([0,1],\Om^0(P_{c_k}))$ such that $f_k(0)=P_{c_k}$ and $f_k(1)=A$ in Theorem 2.6.
}

\medskip

{\bf Proof.} Under the assumption (TCG), it follows from Theorem 1.1 that both $c_1$ and $c_2$ are irrationally elliptic. Then there holds $P_{c_k}\approx R(\tilde{\theta}_{k,1})\,\dm\,R(\tilde{\theta}_{k,2})\,\dm\,R(\tilde{\theta}_{k,3})$ for some $\frac{\tilde{\theta}_{k,1}}{2\pi}, \frac{\tilde{\theta}_{k,2}}{2\pi}, \frac{\tilde{\theta}_{k,3}}{2\pi}\in (0,1)\bs\Q $ for $k=1,2$, respectively.

Then  by (\ref{2.18}), we have
\bea
i(c_k^m)&=&m\left(i(c_k)-3\right)+2\sum_{j=1}^{3}E\left(\frac{m\tilde{\theta}_{k,j}}{2\pi}\right)-3,\quad \nu(c_k^m)=0,\quad k=1,2,\lb{4.2}\\
\hat{i}(c_k)&=&i(c_k)-3+\sum_{j=1}^{3}\frac{\tilde{\theta}_{k,j}}{\pi},\quad k=1,2,\lb{4.3}
\eea
and then by (\ref{2.7}), we get
\bea
n(c_k)=1,\quad k=1,2.\lb{4.4}
\eea
By (\ref{2.6}), (iii) of Lemma 2.2, (\ref{4.2}) and (\ref{4.4}), for the average Euler numbers of $c_1$ and $c_2$ we have
\bea
\hat{\chi}(c_1)=-1,\quad\hat{\chi}(c_2)=-1.\lb{4.5}
\eea

Noticing that $\nu(c_k^m)=0,\,\forall\  m\in\N$ and $k=1,2$, so by Lemma 2.1, (\ref{2.4}), and (iii) of Lemma 2.2, we have
\bea
M_q(k,m)=\left\{\matrix{
   1,&&\quad {\it if}\quad q=i(c_k^m),  \cr
   0,&&\quad {\it if}\quad q\neq i(c_k^m),  \cr}\right.\quad \mbox{for}\  k=1,2,\  m\in\N.\lb{4.6}
\eea

By Claim 3 of the proof of Theorem 1.1, (\ref{4.2}) and (\ref{4.6}), we have
\bea
i(c_2^{2m_2})+\nu(c_2^{2m_2})=i(c_2^{2m_2})=2N-3,\quad M_{2N-3}(2,2m_2)=1.\lb{4.7}
\eea
By (\ref{3.38.0}) and (\ref{4.2}), we have
\bea
i(c_1^{2m_1})+\nu(c_1^{2m_1})=i(c_1^{2m_1})=2N+3,\quad M_{2N+3}(1,2m_1)=1.\lb{4.8}
\eea

It follows from (\ref{4.6})-(\ref{4.8}) and Lemma 3.4 that
\bea
M_q=\sum_{1\le k\le 3}M_q(k,2m_k)=M_q(3,2m_3)\quad \mbox{for}\ 2N-2\le q\le 2N+2\lb{4.9}
\eea
and
\bea
M_{2N+3}(2,2m_2)=0.\lb{4.10}
\eea
 By (\ref{2.8}), Theorem 2.5 and (\ref{4.9}), we have $M_{2N-1}=M_{2N-1}(3,2m_3)\ge b_{2N-1}=1$ and $M_{2N+1}=M_{2N+1}(3,2m_3)\ge b_{2N+1}=1$, then by (\ref{3.16.0}), (\ref{3.17.0}) and Lemma 3.5, there holds
\bea
M_{2N-3}(3,2m_3)=M_{2N+3}(3,2m_3)=0.\lb{4.11}
\eea
By (\ref{4.8}), (\ref{4.10}), (\ref{4.11}) and Lemma 3.4, we have
\bea
M_{2N+3}=\sum_{1 \le k\le 3}M_{2N+3}(k,2m_k+1)+1.\lb{4.12}
\eea
By (\ref{3.18.0}) and Lemma 3.5, $M_{2N+3}(k,2m_k+1)\le 1$ for $k=1,2,3$. Then it follows from (\ref{4.12}) that $M_{2N+3}\le 4$.
On the other hand, by (\ref{2.8}) and Theorem 2.5, we have $M_{2N+3}\ge b_{2N+3}=2$. Thus $M_{2N+3}\in\{2,3,4\}$.

We continue the proof by distinguishing three cases.

\medskip

{\bf Case 1:} $M_{2N+3}=4$.

\medskip

In this case, by (\ref{4.12}), it yields
\bea
M_{2N+3}(k,2m_k+1)=1,\quad\forall\ k=1,2,3.\lb{4.13}
\eea
Then by (\ref{3.18.0}) and Lemma 3.5, we have
\bea
M_{2N+5}(k,2m_k+1)=0,\quad\forall\ k=1,2,3.\lb{4.14}
\eea
It follows from (\ref{4.14}) and Lemma 3.4 that $M_{2N+5}=0$. However by (\ref{2.8}) and Theorem 2.5, we have $M_{2N+5}\ge b_{2N+5}=1$, which is a contradiction.

\medskip

{\bf Case 2:} $M_{2N+3}=3$.

\medskip

In this case, by (\ref{4.12}), one of the following cases may happen.

(i) $M_{2N+3}(3,2m_3+1)=0$ and $M_{2N+3}(1,2m_1+1)=M_{2N+3}(2,2m_2+1)=1$.

(ii) $M_{2N+3}(3,2m_3+1)=1$ and $M_{2N+3}(1,2m_1+1)+M_{2N+3}(2,2m_2+1)=1$.

For (i), by (\ref{4.6}), we have $i(c_1^{2m_1+1})=2N+3$ and $i(c_2^{2m_2+1})=2N+3$, then by (\ref{3.9.0}), it yields $i(c_1)=i(c_2)=3$. So, by (\ref{3.4.0}) and (\ref{4.3}), we get
\bea
5<\hat{i}(c_1)< 6,\quad 5<\hat{i}(c_2)< 6.\lb{4.15}
\eea
By (\ref{3.18.0}), Lemma 3.5 and (i) of Case 2, there holds
\bea
M_q(k,2m_k+1)=0\quad \mbox{for}\ q\neq 2N+3\ \mbox{and}\ k=1,2.\lb{4.16}
\eea

So, according to Lemma 3.4 and (\ref{4.16}),
we have
\bea
M_{2N+5}&=&\sum_{1\le k\le 3}M_{2N+5}(k,2m_k+1)=M_{2N+5}(3,2m_3+1),\lb{4.17}\\
M_{2N+7}&=&\sum_{1\le k\le 3}M_{2N+7}(k,2m_k+1)=M_{2N+7}(3,2m_3+1).\lb{4.18}
\eea
By (\ref{2.8}) and Theorem 2.5, we have $M_{2N+5} \ge b_{2N+5}=1$ and $M_{2N+7} \ge b_{2N+7}=1$, then it follows from (\ref{4.17}) and (\ref{4.18}) that
\bea
M_{2N+5}(3,2m_3+1)\ge 1,\quad M_{2N+7}(3,2m_3+1)\ge 1.
\eea
Thus the assupmtion with $q=6$ in Lemma 3.6 is satisfied, and then by (\ref{3.52.0}) we have
\bea
i(c_3^m)=6m-p_{3,-}-p_{3,0}.\lb{4.20}
\eea
Then we have $n(c_3)=1$ in either case of $P_{c_3}$ by (\ref{2.7}).
By (\ref{4.20}), we have
\bea
\hat{i}(c_3)=6.\lb{4.21}
\eea
By (\ref{2.5}), we have the following  identity
\bea
\sum_{k=1}^{3}\frac{\hat{\chi}(c_k)}{\hat{i}(c_k)}=B(4,1)=-\frac{2}{3}.\lb{4.22}
\eea
Combining (\ref{4.5}), (\ref{4.15}), (\ref{4.21}) and (\ref{4.22}), we obtain $-2<\hat{\chi}(c_3)=-4+\frac{6}{\hat{i}(c_1)}+\frac{6}{\hat{i}(c_2)}<-\frac{8}{5}$, which contradicts to $\hat{\chi}(c_3)\in \Z$, where the latter is due to $n(c_3)=1$ and the definition of $\hat{\chi}(c_3)$.

For (ii), without loss of generality, we assume that $M_{2N+3}(1,2m_1+1)=0$ and $M_{2N+3}(2,2m_1+1)=1$. Thus, similarly, by (\ref{3.18.0}) and  Lemma 3.5 and Lemma 3.4,
we have
\bea
M_{2N+5}=M_{2N+5}(1,2m_1+1),\quad M_{2N+7}=M_{2N+7}(1,2m_1+1).\lb{4.23}
\eea

By (\ref{2.8}) and Theorem 2.5, we have $M_{2N+5} \ge b_{2N+5}=1$ and $M_{2N+7} \ge b_{2N+7}=1$. Then it follows from (\ref{4.23}) that $M_{2N+5}(1,2m_1+1)\ge 1$ and $M_{2N+7}(1,2m_1+1)\ge 1$. Thus by Lemma 3.6, there are five cases for $P_{c_1}$, which contradicts the fact that  $c_1$ is irrationally elliptic.

\medskip

{\bf Case 3:} $M_{2N+3}=2$.

\medskip

In this case, by (\ref{4.12}), one of the following cases may happen:

(i) $M_{2N+3}(3,2m_3+1)=1$ and $M_{2N+3}(1,2m_1+1)=M_{2N+3}(2,2m_2+1)=0$.

(ii) $M_{2N+3}(3,2m_3+1)=0$ and $M_{2N+3}(1,2m_1+1)+M_{2N+3}(2,2m_2+1)=1$.

For (i), by (\ref{3.18.0}) and Lemma 3.5, there holds
\bea
M_q(3,2m_3+1)=0\quad \mbox{for}\  q\neq 2N+3,\lb{4.25}
\eea
and then by Lemma 3.4, for $q=2N+5, 2N+7$, we have
\bea
M_q=\sum_{1\le k \le 3}M_q(k,2m_k+1)=M_q(1,2m_1+1)+M_q(2,2m_2+1).
\eea

On the other hand, by (\ref{2.8}) and Theorem 2.5, we
have $M_{2N+5}\ge b_{2N+5}=1$ and $M_{2N+7}\ge b_{2N+7}=1$, which, together with (\ref{4.6}), implies
\bea
M_q(1,2m_1+1)=\left\{\matrix{
   1,\quad {\it if}\quad q=2N+5,  \cr
   0,\quad {\it if}\quad q\neq 2N+5,  \cr}\right.
\quad M_q(2,2m_2+1)=\left\{\matrix{
   1,\quad {\it if}\quad q=2N+7,  \cr
   0,\quad {\it if}\quad q\neq 2N+7,  \cr}\right.\lb{4.31.0}
\eea
or
\bea
M_q(1,2m_1+1)=\left\{\matrix{
   1,\quad {\it if}\quad q=2N+7,  \cr
   0,\quad {\it if}\quad q\neq 2N+7,  \cr}\right.
\quad M_q(2,2m_2+1)=\left\{\matrix{
   1,\quad {\it if}\quad q=2N+5,  \cr
   0,\quad {\it if}\quad q\neq 2N+5.  \cr}\right.\lb{4.32.0}
\eea

So it follows from Lemma 3.4, (\ref{4.25}) and (\ref{4.31.0})-(\ref{4.32.0}) that
 \bea
M_{2N+9}= \sum_{1\le k\le 3, \atop 1\le m\le 2}M_{2N+9}(k,2m_k+m)= \sum_{1\le k\le 3}M_{2N+9}(k,2m_k+2).\lb{4.27}
\eea

Without loss of generality,  we assume that (\ref{4.31.0}) holds. So we get $i(c_1^{2m_1+1})=2N+5$ and $i(c_2^{2m_2+1})=2N+7$ by (\ref{4.6}). Then by (\ref{3.9.0}), it yields $i(c_1)=5$ and $i(c_2)=7$.

Since $i(c_2)=7$, by (\ref{4.2}), we have the index iteration formula of $c_2$ as follows
\bea
i(c_2^m)=4m-3+2\sum_{j=1}^{3}E\left(\frac{m\tilde{\theta}_{2,j}}{2\pi}\right).\lb{4.29}
\eea
By (\ref{4.29}), it yields $i(c_2^{2})\ge 11$, then by (\ref{3.9.0}), we get $i(c_2^{2m_2+2})\ge 2N+11$. Thus by (\ref{4.6}), we have
\bea
M_{2N+9}(2,2m_2+2)=0.\lb{4.30}
\eea
By (\ref{3.18.0}) and Lemma 3.5, we have
\bea
M_{2N+9}(k,2m_k+2)\le 1,\quad\forall\  k=1,2,3, \lb{4.31}
\eea
which, together with (\ref{4.27}) and (\ref{4.30}), implies $M_{2N+9}\le 2$. However by (\ref{2.8}) and Theorem 2.5, we have $M_{2N+9}\ge b_{2N+9}=2$, which implies $M_{2N+9}=2$.
So we have
\bea
M_{2N+9}(1,2m_1+2)=M_{2N+9}(3,2m_3+2)=1.\lb{4.32}
\eea
It follows from (\ref{4.6}) and (\ref{4.32}) that
\bea
M_q(1,2m_1+2)=M_q(3,2m_3+2)=0\quad \mbox{for}\  q=2N+11,2N+13.\lb{4.33}
\eea
Combining (\ref{4.25}), (\ref{4.31.0}), (\ref{4.33}) and Lemma 3.4, we obtain
\bea
M_q=\sum_{1\le k\le 3, \atop 1\le m\le 2}M_q(k,2m_k+m)=M_q(2,2m_2+2)\quad \mbox{for}\  q=2N+11,2N+13.
\eea
Then by (\ref{4.6}), we obtain that $M_{2N+11}=0$ or $M_{2N+13}=0$, which contradicts to $M_{2N+11}\ge b_{2N+11}=1$ and $M_{2N+13}\ge b_{2N+13}=1$ by (\ref{2.8}) and Theorem 2.5.

For (ii), without loss of generality,  we assume that
\bea
M_{2N+3}(2,2m_2+1)=0,\quad M_{2N+3}(1,2m_1+1)=1.\lb{4.35}
\eea
By (\ref{4.6}), we have $i(c_1^{2m_1+1})=2N+3$ and
\bea
M_q(1,2m_1+1)=0,\quad \forall\ q\neq 2N+3.\lb{4.36}
\eea
Then by (\ref{3.9.0}), it yields $i(c_1)=3$. Then by (\ref{4.2}), we obtain
\bea
i(c_1^m)=2\sum_{j=1}^{3}E\left(\frac{m\tilde{\theta}_{1,j}}{2\pi}\right)-3.\lb{4.37}
\eea
By (\ref{4.37}), it yields $i(c_1^2)\le 9$. Then by (\ref{3.3.0}), we get $i(c_1^2)=9$. And then $i(c_1^{2m_1+2})=2N+9$ by (\ref{3.9.0}). Thus by (\ref{4.6}), we have
\bea
M_q(1,2m_1+2)=\left\{\matrix{
   1,\quad {\it if}\quad q=2N+9,  \cr
   0,\quad {\it if}\quad q\neq 2N+9.  \cr}\right.\lb{4.38}
\eea

\medskip

{\bf Claim 1:} {\it $i(c_2^{2m_2+1})\le 2N+9$, or equivalently, $i(c_2)\le 9$ by (\ref{3.9.0}).}

\medskip

If $i(c_2^{2m_2+1})>2N+9$, by (\ref{3.5.0}) we get $i(c_2^{2m_2+2})\ge i(c_2^{2m_2+1})>2N+9$. Then  by (\ref{4.6}) we know
\bea
M_q(2,2m_2+1)=M_q(2,2m_2+2)=0\quad \mbox{for}\ q=2N+5,2N+7,2N+9. \lb{4.39}
\eea
Combining Lemma 3.4, (\ref{4.36}) and (\ref{4.39}), we have
\bea
M_{2N+5}&=&\sum_{1\le k\le 3}M_{2N+5}(k,2m_k+1)=M_{2N+5}(3,2m_3+1),\lb{4.40}\\
M_{2N+7}&=&\sum_{1\le k\le 3}M_{2N+7}(k,2m_k+1)=M_{2N+7}(3,2m_3+1).\lb{4.41}
\eea

On the other hand, by (\ref{2.8}) and Theorem 2.5, we have $M_{2N+5}\ge b_{2N+5}=1$  and $M_{2N+7}\ge b_{2N+7}=1$. Then by (\ref{4.40}) and (\ref{4.41}), we know that
\bea
M_{2N+5}(3,2m_3+1)\ge 1,\quad M_{2N+7}(3,2m_3+1)\ge 1.\lb{4.42}
\eea

Thus the assumption with $q=6$ in Lemma 3.6 is satisfied, and then we have the index iteration formula of $c_3$ as follows
\bea
i(c_3^m)=6m-p_{3,-}-p_{3,0}.\lb{4.43}
\eea
Then by the fact $\nu(c_3)=p_{3,-}+2p_{3,0}+p_{3,+}$ and (\ref{3.2.0}), we get
\bea
i(c_3)+\nu(c_3)=6+p_{3.0}+p_{3,+}\le 9.\lb{4.44}
\eea
Then by (\ref{3.9.0}) and (\ref{2.22}), it yields $i(c_3^{2m_3+1})+\nu(c_3^{2m_3+1})\le 2N+9$. So, according to Lemma 3.5 and (\ref{4.42}), there holds
\bea
M_{2N+9}(3,2m_3+1)=0.\lb{4.45}
\eea
 Similar to (\ref{3.73.0}), we obtain
\bea
M_{2N+9}(3,2m_3+2)=M_{2N+3}(3,2m_3+1)=0.\lb{4.46}
\eea
Combining (\ref{4.36}), (\ref{4.38}), (\ref{4.39}), (\ref{4.45}), (\ref{4.46}) and Lemma 3.4, we obtain that
\bea
M_{2N+9}= \sum_{1\le k\le 3, \atop 1\le m\le 2}M_{2N+9}(k,2m_k+m)=M_{2N+9}(1,2m_1+2)=1,\lb{4.47}
\eea
which gives a contradiction $1=M_{2N+9}\ge b_{2N+9}=2$ by (\ref{2.8}) and Theorem 2.5. This finished the proof of Claim 1.

Note that $i(c_2)\neq 3$ by (\ref{4.35}), (\ref{4.6}) and (\ref{3.9.0}), it yields $i(c_2)\in\{5,7,9\}$ by Claim 1 since $i(c_2)$ is odd. Next we have three subcases according to the value of $i(c_2)$.

\medskip

{\bf Subcase 3.1:} $i(c_2)=5$.

\medskip

In this subcase, by (\ref{4.6}) and (\ref{3.9.0}), we have
\bea
M_q(2,2m_2+1)=\left\{\matrix{
   1,\quad {\it if}\quad q=2N+5,  \cr
   0,\quad {\it if}\quad q\neq 2N+5.  \cr}\right.\lb{4.48}
\eea
 By Lemma 3.4, (\ref{4.36}) and (\ref{4.48}), we obtain that
\bea
M_{2N+7}=\sum_{1\le k\le 3}M_{2N+7}(k,2m_k+1)=M_{2N+7}(3,2m_3+1).\lb{4.49}
\eea
Together with $M_{2N+7}\ge b_{2N+7}=1$ by (\ref{2.8}) and Theorem 2.5, we get
\bea
M_{2N+7}(3,2m_3+1)\ge 1.\lb{4.50}
\eea

\medskip

{\bf Claim 2:} {\it $M_q(3,2m_3+1)=0$ for $q\ge 2N+11$.}

\medskip

If $M_{q_0}(3,2m_3+1)\ge 1$ for some $q_0\ge 2N+11$, then by (\ref{4.50}) and Lemma 3.5, we know that $i(c_3^{2m_3+1})\le 2N+6$ and   $i(c_3^{2m_3+1})+\nu(c_3^{2m_3+1})\ge q_0+1\ge 2N+12$, which, together with the fact $\nu(c_3^{2m_3+1})\le 6$, implies $i(c_3^{2m_3+1})=2N+6$ and $\nu(c_3^{2m_3+1})=6$ and  $q_0$ only can be $2N+11$. Now by (\ref{3.9.0}) and (\ref{2.22}), we have $i(c_3)=6$ and $\nu(c_3)=6$, which implies that $P_{c_3}\approx I_6$ by (\ref{3.1.0}), (\ref{3.2.0}) and the fact $\nu(c_3)=p_{3,-}+2p_{3,0}+p_{3,+}$. Then $i(c_3)$ must be odd by Proposition 2.7 and the symplectic additivity of symplectic paths. This contradicts to $i(c_3)=6$ and completes the proof of Claim 2.

\medskip

In summary, by Lemma 3.4, (\ref{4.36}), (\ref{4.48}), Claim 2 and (\ref{4.38}), there holds
{\small \bea
&&M_q=\sum_{1\le k\le 3}M_q(k,2m_k+1)=M_q(3,2m_3+1)\quad\mbox{for}\ 2N+4\le q \le 2N+8,\ q\neq 2N+5,\lb{4.51}\\
&&M_{2N+5}=\sum_{1\le k\le 3}M_{2N+5}(k,2m_k+1)=1+M_{2N+5}(3,2m_3+1),\lb{4.52}\\
&&M_{2N+9}=\sum_{1\le k\le 3, \atop 1\le m\le 2}M_{2N+9}(k,2m_k+m)=1+M_{2N+9}(3,2m_3+1)+\sum_{2\le k\le 3}M_{2N+9}(k,2m_k+2),\lb{4.53}\\
&&M_{2N+10}=\sum_{1\le k\le 3, \atop 1 \le m\le 2}M_{2N+10}(k,2m_k+m)=M_{2N+10}(3,2m_3+1)+\sum_{2\le k\le 3}M_{2N+10}(k,2m_k+2),\lb{4.54}\\
&&M_q=\sum_{1\le k\le 3, \atop 1\le m\le 2}M_q(k,2m_k+m)=M_q(2,2m_2+2)+M_q(3,2m_3+2),\  2N+11\le q \le 2N+14.\lb{4.55}
\eea}

Note that  by (\ref{4.2}) and $i(c_2)=5$, we have the index iteration formula of $c_2$ as follows
\bea
i(c_2^m)=2m-3+2\sum_{j=1}^{3}E\left(\frac{m\tilde{\theta}_{2,j}}{2\pi}\right).\lb{4.56}
\eea
It follows from (\ref{4.56})  that $i(c_2^2)\le 13$. Then by (\ref{3.3.0}) and the fact that $i(c_2^2)$ is odd, we get $i(c_2^2)\in \{9,11,13\}$.

We continue the proof by distinguishing three values of $i(c_2^2)$.
\medskip

{\bf Subcase 3.1.1:} $i(c_2^2)=9$.

\medskip

In this subcase, by (\ref{3.9.0}), we have $i(c_2^{2m_2+2})=2N+9$, then by (\ref{4.6}), we have
\bea
M_q(2,2m_2+2)=\left\{\matrix{
   1,\quad {\it if}\quad q=2N+9,  \cr
   0,\quad {\it if}\quad q\neq 2N+9.  \cr}\right.\lb{4.57}
\eea
It follows from (\ref{4.55}) and (\ref{4.57}) that
\bea
M_{2N+11}=M_{2N+11}(3,2m_3+2),\quad M_{2N+13}=M_{2N+13}(3,2m_3+2),\lb{4.58}
\eea
  which, together with $M_{2N+11}\ge b_{2N+11}=1$ and  $M_{2N+13}\ge b_{2N+13}=1$ by (\ref{2.8}) and Theorem 2.5, implies
  \bea
M_{2N+11}(3,2m_3+2)\ge 1,\quad M_{2N+13}(3,2m_3+2)\ge 1.\lb{4.59}
\eea
Thus the assumption with $q=12$ in Lemma 4.1 is satisfied, and then we have
\bea
i(c_3^m) &=& 6m-(p_{3,-}+p_{3,0})-\frac{1+(-1)^m}{2}(q_{3,0}+q_{3,+}),\lb{4.60}\\
\nu(c_3^m) &=& p_{3,-}+2p_{3,0}+p_{3,+}+\frac{1+(-1)^m}{2}(q_{3,-}+2q_{3,0}+q_{3,+}).\lb{4.61}
\eea
Then we can know that
\bea
\hat{i}(c_3)=6\lb{4.62}
\eea
and
\bea
n(c_3)\in\{1,2\}.\lb{4.63}
\eea
By (\ref{3.2.0}), we obtain $i(c_3^2)+\nu(c_3^2)=12+p_{3,0}+p_{3,+}+q_{3,-}+q_{3,0}\le 15$, which implies
\bea
i(c_3^{2m_3+2})+\nu(c_3^{2m_3+2})\le 2N+15\lb{4.64}
\eea
by (\ref{3.9.0}) and (\ref{2.22}). Then by Lemma 3.5 and (\ref{4.59}), there holds
\bea
M_{2N+15}(3,2m_3+2)=0,\quad M_{2N+17}(3,2m_3+2)=0. \lb{4.65}
\eea
It follows from Lemma 3.4, (\ref{4.36}), (\ref{4.38}), (\ref{4.48}), Claim 2, (\ref{4.57}) and (\ref{4.65}) that
\bea
M_{2N+15}&=&\sum_{1\le k\le 3, \atop 1\le m\le 3}M_{2N+15}(k,2m_k+m)=\sum_{1\le k\le 3}M_{2N+15}(k,2m_k+3),\lb{4.66}\\
M_{2N+17}&=&\sum_{1\le k\le 3, \atop 1\le m\le 3}M_{2N+17}(k,2m_k+m)=\sum_{1\le k\le 3}M_{2N+17}(k,2m_k+3).\lb{4.67}
\eea
By (\ref{3.18.0}) and Lemma 3.5, we have
\bea
M_{2N+15}(k,2m_k+3)\le 1,\quad\forall\ k=1,2,3.\lb{4.68}
\eea
Then by (\ref{4.66}) we get $M_{2N+15}\le 3$. We claim that $M_{2N+15}\neq 3$. In fact, if $M_{2N+15}=3$, we have $M_{2N+15}(k,2m_k+3)= 1,\, k=1,2,3$. Then by (\ref{3.18.0}) and Lemma 3.5, there holds $M_{2N+17}(k,2m_k+3)=0,\, k=1,2,3$, which, together with (\ref{4.67}), implies $M_{2N+17}=0$. This gives a contradiction $0=M_{2N+17}\ge b_{2N+17}=1$ by
(\ref{2.8}) and Theroem 2.5. Hence $M_{2N+15}\le 2$. However, again by (\ref{2.8}) and Theroem 2.5, we have $M_{2N+15}\ge b_{2N+15}=2$. So we get $M_{2N+15}=2$.

In summary, in Case 3, we have
\bea
M_{2N+3}=b_{2N+3},\quad M_{2N+15}=b_{2N+15}.\lb{4.69}
\eea
Combining (\ref{4.69}), (\ref{2.8}) and Theorem 2.5, we obtain that
\bea
\sum_{q=2N+4}^{2N+14}(-1)^{q}M_{q}=\sum_{q=2N+4}^{2N+14}(-1)^{q}b_{q}=-6.\lb{4.70}
\eea
By (\ref{4.51})-(\ref{4.55}), (ii) of Case 3, Claim 2, (\ref{4.57}) and (\ref{4.65}), we have
\bea
\sum_{q=2N+4}^{2N+14}(-1)^{q}M_{q}&=&-2+\sum_{q=2N+4}^{2N+10}(-1)^{q}M_{q}(3,2m_3+1)+
\sum_{2N+9}^{2N+14}(-1)^{q}\sum_{2\le k\le 3}M_q(k,2m_k+2)\nn\\
&=&-3+\sum_{q=2N+3}^{2N+15}(-1)^{q}M_{q}(3,2m_3+1)+
\sum_{2N+9}^{2N+15}(-1)^{q}M_q(3,2m_3+2).\lb{4.71}
\eea

Note that $2N+3\le i(c_3^{2m_3+1})\le i(c_3^{2m_3+1})+\nu(c_3^{2m_3+1})\le 2N+15$ and $2N+9\le i(c_3^{2m_3+2})\le i(c_3^{2m_3+2})+\nu(c_3^{2m_3+2})\le 2N+15$ by (\ref{3.18.0}), (\ref{3.5.0}) and (\ref{4.64}), then according to Lemma 2.1 and (\ref{2.4}), we obtain
\bea
&&\sum_{q=2N+3}^{2N+15}(-1)^{q}M_{q}(3,2m_3+1)+
\sum_{2N+9}^{2N+15}(-1)^{q}M_q(3,2m_3+2)\nn\\
&=&\sum_{0\le l \le 6}(-1)^{i(c_3^{2m_3+1})+l}k_l^{\ep(c_3^{2m_3+1})}(c_3^{2m_3+1})+\sum_{0\le l \le 6}(-1)^{i(c_3^{2m_3+2})+l}k_l^{\ep(c_3^{2m_3+2})}(c_3^{2m_3+2}).\lb{4.72}
\eea
Combining (\ref{4.70}), (\ref{4.71}) and (\ref{4.72}), by (\ref{2.6}), we get
\bea
\chi(c_3^{2m_3+1})+\chi(c_3^{2m_3+2})&=&\sum_{0\le l \le 6}(-1)^{i(c_3^{2m_3+1})+l}k_l^{\ep(c_3^{2m_3+1})}(c_3^{2m_3+1})\nn\\
&&+\sum_{0\le l \le 6}(-1)^{i(c_3^{2m_3+2})+l}k_l^{\ep(c_3^{2m_3+2})}(c_3^{2m_3+2})=-3.\lb{4.73}
\eea
Note that by (\ref{4.63}) and (iv) of Lemma 2.2, there holds
\bea
k_j^{\ep(c_3^{2m_3+1})}(c_3^{2m_3+1})=k_j^{\ep(c_3)}(c_3),\quad
k_j^{\ep(c_3^{2m_3+2})}(c_3^{2m_3+2})=k_j^{\ep(c_3^{2})}(c_3^{2}),\quad\forall\  0\le j\le 6.\lb{4.74}
\eea
Then by (\ref{2.6}) it yields
\bea
\chi(c_3)=\chi(c_3^{2m_3+1}),\quad\chi(c_3^2)=\chi(c_3^{2m_3+2}),\lb{4.75}
\eea
which, together with (\ref{4.73}) and (\ref{4.63}), implies
\bea
\chi(c_3)\neq\chi(c_3^2)\lb{4.76}
\eea
since $\chi(c_3^m)\in\Z$,
and
\bea
n(c_3)=2.\lb{4.77}
\eea
 It follows from (\ref{2.6}), (\ref{4.73}), (\ref{4.75}) and (\ref{4.77}) that
\bea
\hat{\chi}(c_3)=\frac{1}{2}\left(\chi(c_3)+\chi(c_3^2)\right)
=\frac{1}{2}\left(\chi(c_3^{2m_3+1})+\chi(c_3^{2m_3+2})\right)=-\frac{3}{2}.\lb{4.78}
\eea

Note that $\hat{i}(c_k)>5$, $k=1,2$ by (\ref{3.4.0}), so it follows from (\ref{4.5}), (\ref{4.62}) and (\ref{4.78}) that
\bea
\frac{\hat{\chi}(c_1)}{\hat{i}(c_1)}=-\frac{1}{\hat{i}(c_1)}>-\frac{1}{5},\quad
\frac{\hat{\chi}(c_2)}{\hat{i}(c_2)}=-\frac{1}{\hat{i}(c_2)}>-\frac{1}{5},\quad \frac{\hat{\chi}(c_3)}{\hat{i}(c_3)}=-\frac{1}{4},\lb{4.79}
\eea
which, together with (\ref{2.5}), yields
\bea -\frac{2}{3}=\frac{\hat{\chi}(c_1)}{\hat{i}(c_1)}+\frac{\hat{\chi}(c_2)}{\hat{i}(c_2)}+\frac{\hat{\chi}(c_3)}{\hat{i}(c_3)}>-\frac{13}{20}.\lb{4.80}
\eea
This is a contradiction.
\medskip

{\bf Subcase 3.1.2:} $i(c_2^2)=11$.

\medskip

In this subcase, by (\ref{3.9.0}), it yields $i(c_2^{2m_2+2})=2N+11$ and then by (\ref{4.6}), we have
\bea
M_q(2,2m_2+2)=\left\{\matrix{
   1,\quad {\it if}\quad q=2N+11,  \cr
   0,\quad {\it if}\quad q\neq 2N+11,  \cr}\right.\lb{4.81}
\eea
which, together with (\ref{4.55}), implies
\bea
M_{2N+13}=M_{2N+13}(3,2m_3+2).\lb{4.82}
\eea
However by (\ref{2.8}) and Theorem 2.5, we have $M_{2N+13}\ge b_{2N+13}=1$,
so we have
\bea
M_{2N+13}(3,2m_3+2)\ge 1.\lb{4.83}
\eea
 Thus by (\ref{3.18.0}) and Lemma 3.5, there holds
 \bea
M_{2N+9}(3,2m_3+2)=0.\lb{4.84}
\eea

Comnining (\ref{4.53}), (\ref{4.81}) and (\ref{4.84}), we obtain that
\bea
M_{2N+9}=1+M_{2N+9}(3,2m_3+1).\lb{4.85}
\eea
Since $M_{2N+9}\ge b_{2N+9}=2$ by (\ref{2.8}) and Theorem 2.5, by (\ref{4.85}), we know that
\bea
M_{2N+9}(3,2m_3+1)\ge 1.\lb{4.86}
\eea
 Noticing that we also have $M_{2N+7}(3,2m_3+1)\ge 1$ by (\ref{4.50}), then by Lemma 3.6, we obtain that
\bea
i(c_3^m)=8m-p_{3,-}-p_{3,0}.\lb{4.87}
\eea
and
\bea
n(c_3)=1.\lb{4.88}
\eea
Similar to (\ref{3.73.0}), we have
\bea
M_{2N+5}(3,2m_3+1)&=&k_{2N+5-i(c_3^{2m_3+1})}^{\ep(c_3^{2m_3+1})}(c_3^{2m_3+1})\nn\\
&=&k_{2N+5-i(c_3^{2m_3+1})}^{\ep(c_3^{2m_3+2})}(c_3^{2m_3+2})\nn\\
&=&M_{2N+5-i(c_3^{2m_3+1})+i(c_3^{2m_3+2})}(3,2m_3+2)\nn\\
&=&M_{2N+13}(3,2m_3+2).\lb{4.89}
\eea

On the other hand, by (\ref{4.87}) and (\ref{3.2.0}), it yields $i(c_3)\ge 5$, then by (\ref{3.9.0}), we have $i(c_3^{2m_3+1})\ge 2N+5$. So, by (\ref{4.86}) and Lemma 3.5, there holds
\bea
M_{2N+5}(3,2m_3+1)=0,\lb{4.90}
\eea
which, together with (\ref{4.89}), contradicts to (\ref{4.83}).
\medskip

{\bf Subcase 3.1.3:} $i(c_2^2)=13$.

\medskip

In this subcase, by (\ref{3.9.0}), it yields $i(c_2^{2m_2+2})=2N+13$ and then by (\ref{4.6}), we have
\bea
M_q(2,2m_2+2)=\left\{\matrix{
   1,\quad {\it if}\quad q=2N+13,  \cr
   0,\quad {\it if}\quad q\neq 2N+13,  \cr}\right.\lb{4.91}
\eea
which, together with (\ref{4.55}), implies
\bea
M_{2N+11}=M_{2N+11}(3,2m_3+2).\lb{4.92}
\eea
However by (\ref{2.8}) and Theorem 2.5, we have $M_{2N+11}\ge b_{2N+11}=1$,
so we have
\bea
M_{2N+11}(3,2m_3+2)\ge 1.\lb{4.93}
\eea
 Thus by (\ref{3.18.0}) and Lemma 3.5, there holds
 \bea
M_{2N+9}(3,2m_3+2)=0.\lb{4.94}
\eea

Similar to Subcase 3.1.2, we have
\bea
i(c_3^m)=8m-p_{3,-}-p_{3,0},\lb{4.95}
\eea
and we obtain a contradiction
\bea
0=M_{2N+3}(3,2m_3+1)=M_{2N+11}(3,2m_3+2)\ge 1.\lb{4.96}
\eea

\medskip

{\bf Subcase 3.2:} $i(c_2)=7$.

\medskip

In this subcase, by (\ref{3.9.0}), it yields $i(c_2^{2m_2+1})=2N+7$ and then by (\ref{4.6}), we have
\bea
M_q(2,2m_2+1)=\left\{\matrix{
   1,\quad {\it if}\quad q=2N+7,  \cr
   0,\quad {\it if}\quad q\neq 2N+7,  \cr}\right.\lb{4.97}
\eea
which, together with (\ref{4.36}) and Lemma 3.4, implies
\bea
M_{2N+5}= \sum_{1\le k\le 3}M_{2N+5}(k,2m_k+1)=M_{2N+5}(3,2m_3+1).\lb{4.98}
\eea
However by (\ref{2.8}) and Theorem 2.5, we have $M_{2N+5}\ge b_{2N+5}=1$,
so we have
\bea
M_{2N+5}(3,2m_3+1)\ge 1.\lb{4.99}
\eea

 \medskip

{\bf Claim 3:} {\it $M_q(3,2m_3+1)=0$, $\forall\ q\ge 2N+9$.}

\medskip

If $M_{q_0}(3,2m_3+1)\ge 1$ for some $q_0\ge 2N+9$, then by (\ref{4.99}) and Lemma 3.5, we know that $i(c_3^{2m_3+1})\le 2N+4$ and   $i(c_3^{2m_3+1})+\nu(c_3^{2m_3+1})\ge q_0+1\ge 2N+10$, which, together with the fact $\nu(c_3^{2m_3+1})\le 6$, implies $i(c_3^{2m_3+1})=2N+4$ and $\nu(c_3^{2m_3+1})=6$ and  $q_0$ only can be $2N+9$. Now by (\ref{3.9.0}) and (\ref{2.22}), we have $i(c_3)=4$ and $\nu(c_3)=6$, which implies that $P_{c_3}\approx I_6$ by (\ref{3.1.0}), (\ref{3.2.0}) and the fact $\nu(c_3)=p_{3,-}+2p_{3,0}+p_{3,+}$. So $i(c_3)$ must be odd by Proposition 2.7 and the symplectic additivity of symplectic paths. This contradicts to $i(c_3)=4$ and completes the proof of Claim 3.

Since $i(c_2)=7$ in this subcase, similar to the proof of (\ref{4.30}), we have
\bea
M_{2N+9}(2,2m_2+2)=0.\lb{4.100}
\eea
Then by Lemma 3.4, (\ref{4.36}), (\ref{4.97}), Claim 3, (\ref{4.38}) and (\ref{4.100}), we obtain
\bea
M_{2N+9}=\sum_{1\le k\le 3, \atop 1\le m\le 2}M_{2N+9}(k,2m_k+m)=1+M_{2N+9}(3,2m_3+2),\lb{4.101}
\eea
which, together with $M_{2N+9}\ge b_{2N+9}=2$ by (\ref{2.8}) and Theorem 2.5, implies that
\bea
M_{2N+9}(3,2m_3+2)\ge 1.\lb{4.102}
\eea
Then by (\ref{3.18.0}) and Lemma 3.5 we get
\bea
M_q(3,2m_3+2)=0,\quad \forall\ q\neq 2N+9. \lb{4.103}
\eea

By Lemma 3.4, (\ref{4.36}), (\ref{4.97}), Claim 3, (\ref{4.38}) and (\ref{4.103}), we obtain
\bea
M_q=\sum_{1\le k\le 3, \atop 1\le m\le 2}M_q(k,2m_k+m)=M_q(2,2m_2+2), \quad\mbox{for}\  q=2N+11, 2N+13.\lb{4.104}
\eea
So there holds $M_{2N+11}=0$ or $M_{2N+13}=0$ by (\ref{4.6}), which contradicts to $M_{2N+11}\ge b_{2N+11}=1$ and $M_{2N+13}\ge b_{2N+13}=1$ by (\ref{2.8}) and Theorem 2.5.

\medskip

{\bf Subcase 3.3:} $i(c_2)=9$.

\medskip

In this subcase, by (\ref{3.9.0}), it yields $i(c_2^{2m_2+1})=2N+9$ and then by (\ref{4.6}), we have
\bea
M_q(2,2m_2+1)=\left\{\matrix{
   1,\quad {\it if}\quad q=2N+9,  \cr
   0,\quad {\it if}\quad q\neq 2N+9,  \cr}\right.\lb{4.105}
\eea
which, together with (\ref{4.36}) and Lemma 3.4, implies
\bea
M_q=\sum_{1\le k\le 3}M_q(k,2m_k+1)=M_q(3,2m_3+1),\quad\mbox{for}\  q=2N+5, 2N+7.\lb{4.106}
\eea

However by (\ref{2.8}) and Theorem 2.5, we have $M_{2N+5}\ge b_{2N+5}=1$ and $M_{2N+7}\ge b_{2N+7}=1$,
so we have $M_{2N+5}(3,2m_3+1)\ge 1$ and $M_{2N+7}(3,2m_3+1)\ge 1$. Thus the assumption with $q=6$ of Lemma 3.6 is satisfied, and then by Lemma 3.6 we conclude that $c_3$ must belong to one of the classes in Theorem 4.2 and the index of $c_3$ is the following
\bea
i(c_3)=6-p_{3,-}-p_{3,0}.\lb{4.107}
\eea

This completes the proof of Theorem 4.2. \hfill\hb

\medskip

\noindent {\bf Acknowledgments}

\medskip

The authors would like to thank sincerely Professor Yiming Long and Professor Wei Wang for their careful reading of the manuscript and valuable comments.

\bibliographystyle{abbrv}

\begin{thebibliography}{100}
\addtolength{\itemsep}{-0.3ex}

\bibitem[Ano74]{Ano} D. V. Anosov, Geodesics in Finsler geometry. {\it Proc.
I.C.M.} (Vancouver, B.C. 1974), Vol. 2. 293-297 Montreal (1975)
(Russian), {\it Amer. Math. Soc. Transl.} 109 (1977) 81-85.
\bibitem[BTZ82]{BTZ1} W. Ballmann, G. Thobergsson and W. Ziller, Closed geodesics
on positively curved manifolds. {\it Ann. of Math.} 116 (1982), 213-247.
\bibitem[BTZ83]{BTZ2} W. Ballmann, G. Thobergsson and W. Ziller, Existence of closed
geodesics on positively curved manifolds. {\it J. Diff. Geom.} 18 (1983), 221-252.
\bibitem[Ban93]{Ban} V. Bangert, On the existence of closed geodesics on two-spheres.
{\it Inter. J. Math.} 4 (1993), 1-10.
\bibitem[BL10]{BaL} V. Bangert and Y. Long,   The existence of two closed
geodesics on every Finsler 2-sphere. {\it Math. Ann.}  346 (2010), 335-366.
\bibitem[Bot56]{Bot} R. Bott, On the iteration of closed geodesics and the Sturm
intersection theory. {\it Comm. Pure Appl. Math.}  9 (1956), 171-206.
\bibitem[Cha93]{Cha} K. C. Chang, Infinite Dimensional Morse Theory and
Multiple Solution Problems. Birkh\"auser. Boston. 1993.
\bibitem[Dua15]{Dua1} H. Duan, Non-hyperbolic closed geodesics on positively curved Finsler spheres. {\it J. Funct. Anal.}  269 (2015), 3645-3662.
\bibitem[Dua16]{Dua2} H. Duan, Two elliptic closed geodesics on positively curved Finsler spheres. {\it J. Diff. Equa.} 260 (2016), 8388-8402.
\bibitem[DL16]{DL16} H. Duan and H. Liu, Closed geodesics on positively curved Finsler 3-spheres. {\it Adv. Nonlinear Stud.} 16 (2016), 159-171.
\bibitem[DuL10]{DuL} H. Duan and Y. Long, The index growth and mutiplicity
of closed geodesics. {\it J. Funct. Anal.} 259 (2010), 1850-1913.
\bibitem[DLW16]{DLW} H. Duan, Y. Long and W. Wang,
The enhanced common index jump theorem for symplectic paths and non-hyperbolic closed geodesics on Finsler manifolds.  {\it Calc. Var. and PDEs.} 55 (2016), Art. 145, 28 pp.
\bibitem[Hin84]{Hin1} N. Hingston,  Equivariant Morse theory and closed geodesics.
{ \it J. Diff. Geom.} 19 (1984), 85-116.
\bibitem[Hin93]{Hin2} N. Hingston,  On the growth of the number of closed geodesics
on the two-sphere. {\it Inter. Math. Research Notices.} 9 (1993), 253-262.
\bibitem[Hin97]{Hin3} N. Hingston,  On the length of closed geodesics on a two-sphere.
{\it Proc. Amer. Math. Soc.} 125 (1997), 3099-3106.
\bibitem[Fra92]{Fra} J. Franks,  Geodesics on $S^2$ and periodic points of
annulus diffeomorphisms. {\it Invent. Math.}  108 (1992), 403-418.
\bibitem[Kat73]{Kat} A. B. Katok,  Ergodic properties of degenerate integrable
Hamiltonian systems. {\it Izv. Akad. Nauk SSSR.} 37 (1973) (Russian), {\it Math. USSR-Isv.}
 7 (1973), 535-571.
\bibitem[Liu05]{Liu} C. Liu,  The relation of the Morse index of closed geodesics with
the Maslov-type index of symplectic paths. {\it Acta Math. Sin. (Engl. Ser.)} 21 (2005), 237-248.
\bibitem[Lon99]{Lon1} Y. Long,  Bott formula of the Maslov-type index theory.
{\it Pacific J. Math.} 187 (1999), 113-149.
\bibitem[Lon00]{Lon2} Y. Long,  Precise iteration formulae of the Maslov-type index
theory and ellipticity of closed characteristics.  {\it Adv. Math.} 154 (2000), 76-131.
\bibitem[Lon02]{Lon3} Y. Long,  Index Theory for Symplectic Paths with Applications.
Progress in Math. 207, Birkh\"auser. 2002.
\bibitem[Lon06]{Lon4} Y. Long, Multiplicity and stability of closed geodesics on Finsler
2-spheres. {\it J. Euro. Math. Soc.} 8 (2006), 341-353.
\bibitem[LD09]{LoD} Y. Long and H. Duan,  Multiple closed geodesics on 3-spheres.
{\it Adv. Math.} 221 (2009), 1757-1803.
\bibitem[LW08]{LW}Y. Long and W. Wang, Stability of closed geodesics on Finsler 2-spheres. {\it J. Funct. Anal.} 255 (2008), 620-641.
\bibitem[LZ02]{LoZ} Y. Long and C. Zhu, Closed charateristics on compact convex hypersurfaces
in $\R^{2n}$. {\it Ann. of Math.} 155 (2002), 317-368.
\bibitem[Rad89]{Rad} H.-B. Rademacher,  On the average indices of closed
geodesics. {\it J. Diff. Geom.} 29 (1989), 65-83.
\bibitem[Rad92]{Rad2} H.-B. Rademacher, Morse Theorie und geschlossene
Geodatische. {\it Bonner Math. Schriften} Nr. 229 (1992).
\bibitem[Rad04]{Rad3} H.-B. Rademacher, A sphere theorem for non-reversible Finsler metric. {\it Math. Ann.} 328 (2004), 373-387.
\bibitem[Rad07]{Rad4} H.-B. Rademacher, Existence of closed geodesics on positively curved
Finsler manifolds.  {\it Ergod. Th. \& Dynam. Sys.} 27 (2007), 957-969.
\bibitem[Wan12]{Wan} W. Wang, On a conjecture of Anosov. {\it Adv. Math.} 230 (2012), 1597-1617.
\bibitem[Zil82]{Zil} W. Ziller,  Geometry of the Katok examples. {\it Ergod. Th. \& Dynam. Sys.}
3 (1982), 135-157.

\end{thebibliography}

\end{document}